\documentclass[english]{amsart}

\usepackage{amsmath}

\usepackage{amssymb}

\usepackage{amscd} \usepackage{tabularx}
\usepackage[all,cmtip,line]{xy}

\newcommand{\ra}{\rightarrow}
\newcommand{\lra}{\longrightarrow}

\newcommand{\into}{\hookrightarrow}

\newcommand{\da}{\downarrow}
\newcommand{\iso}{\stackrel{\sim}{\ra}}
\newcommand{\liso}{\stackrel{\sim}{\lra}}

\newlength{\ownl}

\newcommand{\norm}{{\mbox{\bf N}}}
\newcommand{\ndiv}{{\mbox{$\not| $}}}

\newcommand{\Art}{{\operatorname{Art}\,}}
\newcommand{\Aut}{{\operatorname{Aut}\,}}
\newcommand{\BC}{{\operatorname{BC}\,}}

\newcommand{\End}{{\operatorname{End}\,}}

\newcommand{\Frob}{{\operatorname{Frob}}}
\newcommand{\Gal}{{\operatorname{Gal}\,}}
\newcommand{\Groth}{{\operatorname{Groth}\,}}
\newcommand{\Hom}{{\operatorname{Hom}\,}}

\newcommand{\Irr}{{\operatorname{Irr}}}

\newcommand{\Lie}{{\operatorname{Lie}\,}}

\newcommand{\Red}{{\operatorname{Red}}}

\newcommand{\Iw}{{\operatorname{Iw}}}

\newcommand{\WD}{{\operatorname{WD}}}

\newcommand{\Spec}{{\operatorname{Spec}\,}}

\newcommand{\Spp}{{\operatorname{Sp}\,}}

\newcommand{\rec}{{\operatorname{rec}}}
\newcommand{\tr}{{\operatorname{tr}\,}}

\newcommand{\nind}{{\operatorname{n-Ind}\,}}

\newcommand{\st}{{\operatorname{st}}}

\newcommand{\ab}{{\operatorname{ab}}}

\newcommand{\AV}{{\operatorname{AV}}}

\newcommand{\et}{{\operatorname{et}}}

\newcommand{\op}{{\operatorname{op}}}

\newcommand{\semis}{{\operatorname{ss}}}
\newcommand{\Fsemis}{{\operatorname{F-ss}}}

\newcommand{\A}{{\mathbb{A}}}

\newcommand{\C}{{\mathbb{C}}}

\newcommand{\F}{{\mathbb{F}}}
\newcommand{\G}{{\mathbb{G}}}

\newcommand{\K}{{\mathbb{K}}}

\newcommand{\Q}{{\mathbb{Q}}}
\newcommand{\R}{{\mathbb{R}}}

\newcommand{\Z}{{\mathbb{Z}}}

\newcommand{\CA}{{\mathcal{A}}}

\newcommand{\CC}{{\mathcal{C}}}

\newcommand{\CG}{{\mathcal{G}}}
\newcommand{\CH}{{\mathcal{H}}}

\newcommand{\CL}{{\mathcal{L}}}

\newcommand{\CO}{{\mathcal{O}}}

\newcommand{\CR}{{\mathcal{R}}}

\newcommand{\gX}{{\mathfrak{X}}}

\newcommand{\gm}{{\mathfrak{m}}}
\newcommand{\gn}{{\mathfrak{n}}}

\newcommand{\barF}{\overline{{F}}}

\newcommand{\barK}{\overline{{K}}}

\newcommand{\barM}{\overline{{M}}}

\newcommand{\barX}{\overline{{X}}}

\newcommand{\barAA}{\overline{{\A}}}

\newcommand{\barFF}{\overline{{\F}}}

\newcommand{\barQQ}{\overline{{\Q}}}

\newcommand{\tR}{\widetilde{{R}}}

\newcommand{\bareta         }{\overline{\eta}}

\newcommand{\tPi   }{\widetilde{\Pi}}

\newcommand{\proj}{{\operatorname{proj}}}

\newcommand{\Qbar}{{\overline{\Q}}}

\def\RCS$#1: #2 ${\expandafter\def\csname RCS#1\endcsname{#2}}
\RCS$Revision: 54 $
\RCS$Date: 2011-03-23 14:00:36 -0400 (Wed, 23 Mar 2011) $

\newcommand{\To}{\longrightarrow}

\newcommand{\onto}{\twoheadrightarrow}
\newcommand{\isoto}{\stackrel{\sim}{\To}} 
 
\newcommand{\cL}{\mathcal{L}}

\newcommand{\bb}{\mathbb}

\newcommand{\seq}[1]{\left<#1\right>}

\newcommand{\Sp}{\operatorname{Sp}}

\newcommand{\HT}{\operatorname{HT}}

\newcommand{\cO}{\mathcal{O}}

\newcommand{\Qpbar}{\overline{\Q}_p}
\newcommand{\Qlbar}{\overline{\Q}_{l}}

\DeclareMathOperator{\Ram}{Ram}
\DeclareMathOperator{\Spl}{Spl}
\DeclareMathOperator{\Unr}{Unr}
\DeclareMathOperator{\Sh}{Sh}
\DeclareMathOperator{\ur}{ur}
\DeclareMathOperator{\fin}{fin}
\DeclareMathOperator{\Ig}{Ig}
\DeclareMathOperator{\nRed}{n-Red}
\DeclareMathOperator{\LJ}{LJ}
\DeclareMathOperator{\id}{id}

\DeclareMathOperator{\Ma}{Ma}

\usepackage{amsthm}

 \newtheorem{ithm}{Theorem}
\newtheorem{thm}{Theorem}[section]
\newtheorem{corollary}[thm]{Corollary}
\newtheorem{cor}[thm]{Corollary}
 
\newtheorem{lem}[thm]{Lemma} \newtheorem{prop}[thm]{Proposition}
 \theoremstyle{definition}
 \theoremstyle{definition}
 \theoremstyle{remark}

\numberwithin{equation}{subsection}

\theoremstyle{definition}

\begin{document}
\title{Local-global compatibility for $l=p$, I.}

\author{Thomas Barnet-Lamb}\email{tbl@brandeis.edu}\address{Department of Mathematics, Brandeis University}
\author{Toby Gee} \email{gee@math.northwestern.edu} \address{Department of
  Mathematics, Northwestern University} \author{David Geraghty}
\email{geraghty@math.princeton.edu}\address{Princeton University and
  Institute for Advanced Study} \author{Richard Taylor}
\email{rtaylor@math.harvard.edu}\address{Department of Mathematics,
  Harvard University} \thanks{The second author was partially supported
  by NSF grant DMS-0841491, the third author was partially supported
  by NSF grant DMS-0635607 and the fourth author was partially
  supported by NSF grant DMS-0600716 and by the Oswald Veblen and Simonyi Funds at the IAS}  

\begin{abstract}We prove the compatibility of the local and global
  Langlands correspondences at places dividing $l$ for the $l$-adic
  Galois representations associated to regular algebraic conjugate
  self-dual cuspidal automorphic representations of $GL_n$ over an imaginary CM
  field, under the assumption that the automorphic representations
  have Iwahori-fixed vectors at places dividing $l$ and have
  Shin-regular weight.

\end{abstract}
\maketitle
\newpage

\section*{Introduction.}
\label{sec:intro}In this paper we prove the compatibility at places dividing
$l$ of the local and global Langlands correspondences for the $l$-adic
Galois representations associated to regular algebraic 
conjugate self-dual cuspidal automorphic representations of $GL_n$
over an imaginary CM field in the special case that the automorphic representations
  have Iwahori-fixed vectors at places dividing $l$ and have
  Shin-regular weight. In the sequel to this paper
  \cite{blggtlocalglobalII} we build on these results to prove the
  compatibility in general (up to semisimplification in the case of
  non-Shin-regular weight).

Our main result is as follows (see Theorem \ref{thm:semistable Shin-regular case} and Corollary \ref{cor:purity l=p semistable case}).

\begin{ithm}
  \label{thma} Let $m\geq 2$ be an integer, $l$ a rational prime and $\imath:
  \Qlbar \iso \bb C$.  Let $F$ be an imaginary CM field and
  $\Pi$ a regular algebraic, conjugate self-dual cuspidal automorphic representation of $GL_m(\A_F)$. If $\Pi$
  has Shin-regular weight and $v|l$ is a place of $F$ such that
  $\Pi_{v}^{\Iw_{m,v}}\neq \{0\}$, then
\[ \imath\WD(r_{l,\imath}(\Pi)|_{G_{F_v}})^{\Fsemis} \cong \rec(\Pi_{v} \otimes |\det |^{(1-m)/2}). \]
In particular $\WD(r_{l,\imath}(\Pi)|_{G_{F_v}})$ is pure. 
\end{ithm}

(See Section \ref{sec:irreducibility} for any unfamiliar terminology.)
The proof is essentially an immediate application of the methods of
\cite{ty}, applied in the setting of \cite{shin} rather than that of
\cite{ht}, and we refer the reader to the introductions of those
papers for the details of the methods that we use. Indeed, if $\Pi$ is
square-integrable at some finite place, then the result is implicit in
\cite{ty}, although it is not explicitly recorded there. For the
convenience of the reader, we make an effort to make our proof as
self-contained as possible.

\subsection*{Notation and terminology.}

We write all matrix transposes on the left; so ${}^t\!A$ is the
transpose of $A$. We let $B_m\subset GL_m$ denote the Borel subgroup
of upper triangular matrices and $T_m\subset GL_m$ the diagonal
torus. We let $I_m$ denote the identity matrix in $GL_m$. We will
sometimes denote the product $GL_m\times GL_n$ by $GL_{m,n}$.

If $M$ is a field, we let $\barM$ denote an algebraic closure of $M$
and $G_M$ the absolute Galois group $\Gal(\barM/M)$.  Let
$\epsilon_l$ denote the $l$-adic cyclotomic character

Let $p$ be a rational prime and $K/\Q_p$ a finite extension. We let
$\cO_K$ denote the ring of integers of $K$, $\wp_K$ the maximal ideal
of $\cO_K$, $k(\nu_K)$ the residue field $\cO_K/\wp_K$,
$\nu_K:K^\times\onto \Z$ the canonical valuation and
$|\;|_K:K^\times\ra \Q^\times$ the absolute value given by
$|x|_K=\#(k(\nu_K))^{-\nu_K(x)}$. We let $|\;|_K^{1/2}:K^\times \ra
\R^\times_{>0}$ denote the unique positive unramified square root of $|\;|_K$. If $K$ is clear from the context, we
will sometimes write $|\;|$ for $|\;|_K$. We let $\Frob_K$ denote the
geometric Frobenius element of $G_{k(\nu_K)}$ and $I_K$ the kernel of
the natural surjection $G_K\onto G_{k(\nu_K)}$. We will sometimes
abbreviate $\Frob_{\Q_p}$ by $\Frob_p$. We let $W_K$ denote the
preimage of $\Frob_K^\Z$ under the map $G_K\onto G_{k(\nu(K))}$,
endowed with a topology by decreeing that $I_K\subset W_K$ with its
usual topology is an open subgroup of $W_K$. We let $\Art_K : K^\times
\iso W_K^\ab$ denote the local Artin map, normalized to take
uniformizers to lifts of $\Frob_K$.

Let $\Omega$ be an algebraically closed field of characteristic $0$. A
Weil-Deligne representation of $W_K$ over $\Omega$ is a triple
$(V,r,N)$ where $V$ is a finite dimensional vector space over
$\Omega$, $r:W_K\ra GL(V)$ is a representation with open kernel and
$N:V\ra V$ is an endomorphism with
$r(\sigma)Nr(\sigma)^{-1}=|\Art^{-1}_K(\sigma)|_K N$. We say that
$(V,r,N)$ is Frobenius semisimple if $r$ is semisimple and we let
$(V,r,N)^{\Fsemis}$ denote the Frobenius semisimplification
of $(V,r,N)$  (see for instance Section 1 of \cite{ty}) and we let $(V,r,N)^{\semis}$ denote
$(V,r^{\semis},0)$. If $\Omega$ has the
same cardinality as $\C$, we have the notions of a Weil-Deligne
representation being \emph{pure} or \emph{pure of weight $k$} -- see
the paragraph before Lemma 1.4 of \cite{ty}.

We will let $\rec_K$ be the local Langlands correspondence of
\cite{ht}, so that if $\pi$ is an irreducible complex admissible
representation of $GL_n(K)$, then $\rec_K(\pi)$ is a Weil-Deligne
representation of the Weil group $W_K$.  We will write $\rec$ for
$\rec_K$ when the choice of $K$ is clear. If $\rho$ is a continuous
representation of $G_K$ over $\barQQ_l$ with $l\neq p$ then we will
write $\WD(\rho)$ for the corresponding Weil-Deligne representation of
$W_K$. (See for instance Section 1 of \cite{ty}.) 

If $m\geq 1$ is an integer, we let $\Iw_{m,K}\subset GL_m(\cO_K)$
denote the subgroup of matrices which map to an upper triangular
matrix in $GL_m(k(\nu_K))$. If $\pi$ is an irreducible
admissible supercuspidal representation of $GL_m(K)$ and $s\geq 1$ is
an integer we let $\Spp_s(\pi)$ be the square integrable
representation of $GL_{mr}(K)$ defined for instance in Section I.3 of
\cite{ht}. Similarly, if $r : W_K \ra GL_m(\Omega)$ is an irreducible
representation with open kernel and $\pi$ is the supercuspidal
representation $\rec_K^{-1}(r)$, we let $\Sp_s(r)=\rec_K(\Sp_s(\pi))$.
If $K'/K$ is a finite extension and if $\pi$ is an irreducible smooth
representation of $GL_n(K)$ we will write $\BC_{K'/K}(\pi)$ for the
base change of $\pi$ to $K'$ which is characterized by
$\rec_{K'}(\pi_{K'})= \rec_K(\pi)|_{W_{K'}}$.

If $\rho$ is a continuous de Rham representation of $G_K$ over
$\barQQ_p$ then we will write $\WD(\rho)$ for the corresponding
Weil-Deligne representation of $W_K$ (its construction, which is due
to Fontaine, is recalled in Section 1 of \cite{ty}), and if $\tau:K \into \barQQ_p$
is a continuous embedding of fields then we will write
$\HT_\tau(\rho)$ for the multiset of Hodge-Tate numbers of $\rho$ with
respect to $\tau$. Thus $\HT_\tau(\rho)$ is a multiset of $\dim \rho$
integers.  In fact, if $W$ is a de Rham representation of $G_K$ over
$\barQQ_p$ and if $\tau:K \into \barQQ_p$ then the multiset
$\HT_\tau(W)$ contains $i$ with multiplicity $\dim_{\barQQ_p} (W
\otimes_{\tau,K} \widehat{\barK}(i))^{G_K} $. Thus for example
$\HT_\tau(\epsilon_l)=\{ -1\}$.

If $F$ is a number field and $v$ a prime of $F$, we will often denote
$\Frob_{F_v}$, $k(\nu_{F_v})$ and $\Iw_{m,F_v}$ by $\Frob_v$, $k(v)$
and $\Iw_{m,v}$. If $\sigma:F \into \barQQ_p$ or $\C$ is an embedding of fields,
then we will write $F_\sigma$ for the 
closure of the image of $\sigma$. If $F'/F$ is a soluble, finite Galois extension and
if $\pi$ is a cuspidal automorphic representation of $GL_m(\A_F)$ we
will write $\BC_{F'/F}(\pi)$ for its base change to $F'$, an
automorphic representation of $GL_n(\A_{K'})$. If $R:G_F \ra
GL_m(\Qlbar)$ is a continuous representation, we say that $R$ is {\it
  pure of weight $w$} if for all but finitely many primes $v$ of $F$,
$R$ is unramified at $v$ and every eigenvalue of $R(\Frob_v)$ is a
Weil $(\#k(v))^w$-number. (See Section 1 of \cite{ty}.)  If $F$ is an
imaginary CM field, we will denote its maximal totally real subfield
by $F^+$ and let $c$ denote the non-trivial element of $\Gal(F/F^+)$.

\section{Automorphic Galois
  representations}
\label{sec:irreducibility}
\setcounter{subsection}{1}

We recall some now-standard
notation and terminology. 
Let $F$ be an imaginary CM field with maximal totally real subfield $F^+$. 
By a {\em RACSDC} (regular, algebraic, 
conjugate self dual, cuspidal) automorphic representation of
$GL_m(\A_{F})$ we mean that\begin{itemize}
\item $\Pi$ is a cuspidal automorphic representation of $GL_m(\A_F)$ such that  $\Pi_\infty$ has the same infinitesimal character as some irreducible
algebraic representation of the restriction of scalars from $F$ to $\Q$ of
$GL_m$,
\item and $\Pi^c \cong \Pi^\vee$.
\end{itemize}
We will say that
$\Pi$ has {\em level prime to $l$} (resp. {\em level
  potentially prime to $l$}) if for all $v|l$ the representation
$\Pi_v$ is unramified (resp. becomes unramified after a finite base
change).

If $\Omega$ is an algebraically closed field of characteristic $0$ we will
write $(\Z^m)^{\Hom(F,\Omega),+}$ for the set of $a=(a_{\tau,i}) \in (\Z^m)^{\Hom(F,\Omega)}$ satisfying
\[ a_{\tau,1} \geq \dots \geq a_{\tau,m}. \]
We will write $(\Z^m)^{\Hom(F,\Omega)}_0$
for the subset of elements $a \in (\Z^m)^{\Hom(F,\Omega)}$ with 
\[ a_{\tau,i}+a_{\tau \circ c, m+1-i}=0. \]
If $F'/F$ is a finite extension we define $a_{F'} \in (\Z^m)^{\Hom(F',\Omega),+}$ by
\[ (a_{F'})_{\tau,i} = a_{\tau|_F,i}. \] Following \cite{shin} we will
be interested, \emph{inter alia}, in the case that either $m$ is odd;
or that $m$ is even and for some $\tau \in \Hom(F,\Omega)$ and for some
odd integer $i$ we have $a_{\tau,i}>a_{\tau,i+1}$. If either of these
conditions hold then we will say that $a$ is {\em Shin-regular}. (We warn the reader
that this is often referred to as `slightly regular' in the
literature. However as this notion is strictly stronger than
`regularity' we prefer to use the term Shin-regular.) 
If $a \in (\Z^m)^{\Hom(F,\C),+}$, let $\Xi_a$ denote the irreducible algebraic representation of
$GL_m^{\Hom(F,\C)}$ which is the tensor product over $\tau$ of the
irreducible representations of $GL_n$ with highest weights $a_\tau$. We will
say that a RACSDC automorphic representation $\Pi$ of $GL_m(\A_F)$ has
{\em weight $a$} if $\Pi_\infty$ has the same infinitesimal character as
$\Xi_a^\vee$. Note that in this case $a$ must lie in $(\Z^m)^{\Hom(F,\C)}_0$.

We recall (see  Theorem 1.2 of \cite{blght}) that to a RACSDC automorphic representation $\Pi$ of
$GL_m(\A_F)$ and $\imath:\Qlbar \iso \C$ we can associate a continuous semisimple
representation
\[ r_{l,\imath}(\Pi): \Gal(\barF/F) \lra GL_m(\Qlbar) \]
with the properties described in Theorem 1.2 of \cite{blght}. In particular
\[ r_{l,\imath}(\Pi)^c \cong r_{l,\imath}(\Pi)^\vee \otimes \epsilon_l^{1-m}.\] 
For $v|l$ a place of $F$, the representation
$r_{l,\imath}(\Pi)|_{G_{F_v}}$ is de Rham and if $\tau:F \into \barQQ_l$ then
\[ \HT_\tau(r_{l,\imath}(\pi))=\{ a_{\imath\tau,1}+m-1, a_{\imath \tau,2}+m-2,...,a_{\imath \tau,m} \}. \]
If $v\ndiv l$, then the main result of \cite{ana} states that
\[  \imath\WD(r_{l,\imath}(\Pi)|_{G_{F_v}})^{\Fsemis} \cong \rec(\Pi_v \otimes
|\det|^{(1-m)/2}). \]

We recall the following result which will prove useful.

\begin{prop}
  \label{thm:redn to ss}
  Let $\Omega$ be an algebraically closed field of characteristic 0
  and of the same cardinality as $\C$. 
  \begin{enumerate}
  \item Suppose $K/\Q_p$ is a finite extension. Let $(V,r,N)$ and $(V',r',N')$ be pure, Frobenius semisimple
    Weil-Deligne representations of $W_K$ over $\Omega$. If the representations $ (V,r^{\semis})$ and $ (V',(r')^{\semis})
    $ are isomorphic, then $(V,r,N)\cong (V',r',N')$.
  \item If $F$ is an imaginary CM field and $\Pi$ is a RACSDC
    automorphic representation of $GL_n(\A_F)$, then for each
    $\imath:\Omega\iso\C$ and each finite place $v$ of $F$,
    $\imath^{-1}\rec(\Pi_v)$ is pure.
  \end{enumerate}
\end{prop}

\begin{proof}
  The first part follows from Lemma 1.4(4) of \cite{ty}. For the
  second part, Theorem 1.2 of \cite{ana} states that
  $\Pi_v$ is tempered for each finite place $v$ of $F$. If $\sigma$ is
  an automorphism of $\C$, then there is a RACSDC automorphic
  representation $\Pi'={}^\sigma\Pi^\infty \otimes \Pi'_\infty$ of
  $GL_n(\A_F)$ (see Th\'eor\`eme 3.13 of \cite{MR1044819}) and we
  deduce that $\sigma\Pi_v$ is tempered. The second part then follows
  from this and Lemma 1.4(3) of \cite{ty}.
\end{proof}

We can now state our main results.

\begin{thm}
  \label{thm:semistable Shin-regular case}
  Let $m\geq 2$ be an integer, $l$ a rational prime and $\imath:
  \Qlbar \iso \bb C$.  Let $L$ be an imaginary CM field and
  $\Pi$ a RACSDC automorphic representation of $GL_m(\A_L)$. If $\Pi$
  has Shin-regular weight and $v|l$ is a place of $L$ such that
  $\Pi_{v}^{\Iw_{m,v}}\neq \{0\}$, then
\[ \imath\WD(r_{l,\imath}(\Pi)|_{G_{L_v}})^{\Fsemis} \cong \rec(\Pi_{v} \otimes |\det |^{(1-m)/2}). \]
\end{thm}

Before turning to the proof, we first record a corollary.

\begin{cor}
  \label{cor:purity l=p semistable case}
  Let $m\geq 2$ be an integer, $l$ a rational prime and $\imath:
  \Qlbar \iso \bb C$.  Let $L$ be an imaginary CM field and $\Pi$ a RACSDC
  automorphic representation of $GL_m(\A_L)$. If $\Pi$ has
  Shin-regular weight and $v|l$ is a place of $L$ then
  $\WD(r_{l,\imath}(\Pi)|_{G_{L_v}})$ is pure.
\end{cor}

\begin{proof}
Choose a finite CM
  soluble Galois extension
  $F/L$ such that for each prime $w|v$ of $F$,
  $\BC_{F_w/L_v}(\Pi_{v})^{\Iw_{m,w}}\neq\{0\}$. Then
  $\WD(r_{l,\imath}(\Pi)|_{G_{F_w}})$ is pure by Theorem
  \ref{thm:semistable Shin-regular case} and Proposition \ref{thm:redn to ss}.
  Lemma 1.4 of \cite{ty} then implies that
  $\WD(r_{l,\imath}(\Pi)|_{G_{L_v}})$ is pure.
\end{proof}

The rest of this paper will be devoted to the proof of Theorem \ref{thm:semistable Shin-regular case}.

\section{Notation and running assumptions}
\label{subsec:notat-runn-assumpt}

For the convenience of the reader, we recall here
the following notation which appears in \cite{shin}:

\begin{itemize}
\item If $\Pi$ is a RACSDC automorphic representation of $GL_m(\A_F)$
  for some integer $m\geq 2$ and an imaginary CM field $F$, or if $\Pi$ is an
  algebraic Hecke character of $\A_M^\times/M^\times$ for a number
  field $M$, then $R_{l,\imath}(\Pi)$ denotes $r_{l,\imath}(\Pi^\vee)$.
\item If $L/F$ is a finite extension of number fields, then $\Ram_{L/F}$
  (resp.\ $\Unr_{L/F}$, resp.\ $\Spl_{L/F}$) denotes the set of finite
  places of $F$ which are ramified (resp.\ unramified, resp.\
  completely split) in $L$. We denote by $\Spl_{L/F,\Q}$ the set
  of rational primes $p$ such that every place of $F$ above $p$ 
  splits completely in $L$.
\item If $F$ is a number field and $\pi$ is an automorphic
  representation of $GL_{n/F}$, then $\Ram_{\Q}(\pi)$ denotes the set
  of rational primes $p$ such that there exists a place $v|p$ of $F$
  with $\pi_v$ ramified.
\item If $G$ is a group of the form $H(F)$ for
  $F/\Q_p$ finite and $H/F$ a reductive group; or $H(\A_F^T)$ for $F$
  a number field, $H/F$ a reductive group and $T$ a finite set of
  places of $F$ containing all infinite places; or a product of groups
  of this form, then we let $\Irr(G)$ (resp.\ $\Irr_l(G)$) denote the set of
  isomorphism classes of irreducible admissible representations of $G$
  on $\C$-vector spaces (resp.\ $\Qlbar$-vector spaces). We let
  $\Groth(G)$ (resp.\ $\Groth_l(G)$) denote the Grothendieck group of
  the category of admissible $\C$-representations (resp.\ $\Qlbar$-representations) of $G$. (See Section I.2 of \cite{ht}.)
\item $\epsilon : \Z \ra \{0,1\}$ is the unique function such that
  $\epsilon(n)\equiv n \mod 2$.
\item $\Phi_n$ is the matrix in $GL_n$ with
  $(\Phi_n)_{ij}=(-1)^{i+1}\delta_{i,n+1-j}$.
\item If $R\to S$ is a homomorphism of commutative rings, $R_{S/R}$ denotes the
  restriction of scalars functor.
\item If $\pi$ is a representation of a group $G$ with a central
  character, we denote the central character by $\psi_\pi$.
\end{itemize}

We now fix the following notations and assumptions which will be in force from Section
\ref{subsec:shimura-varieties} to Section \ref{subsec:comp-cohom-Y}:
\begin{itemize}
\item $E$ is a quadratic imaginary field;
\item $F^+$ is a totally real field with $[F^+:\Q]\geq 2$;
\item $F=EF^+$ and $\Ram_{F/\Q}\subset \Spl_{F/F^+,\Q}$;
\item $\tau : F\into\C$ is an embedding and $\tau_E=\tau|_E$;
\item $\Phi_\C=\Hom(F,\C)$ and $\Phi^+_\C=\Hom_{E,\tau_E}(F,\C)$;
\item $n\geq 3$ is an odd integer;
\item $p \in \Spl_{E/\Q}$ is a rational prime and $u|p$ is a prime of $E$;
\item $w$ is a prime of $F$ above $u$ and $w_1=w,w_2,\ldots,w_r$
  denote all of the primes of $F$ above $u$;
\item $\iota_p:\Qpbar \iso \C$ is an isomorphism such that
  $\iota_p^{-1}\circ \tau$ induces the place $w$;
\item $l$ is a rational prime (possibly equal to $p$) and $\imath :
  \Qlbar \iso \C$.
\end{itemize}
Define algebraic groups $G_n$ and $\G_n$ over $\Z$ by setting
\begin{align*} G_{ n}(R)&=\{ (\lambda,g) \in R^\times\times
GL_{ n}(\CO_F\otimes_{\Z}R) : g\Phi_{ n}{}^t\!g^c = \lambda
\Phi_{ n}\}, \textrm{ and} \\ \G_{ n}(R)&= R_{\CO_E/\Z}(G_n\times_\Z \CO_E)(R)=G_{ n}(\CO_E \otimes_\Z R)
\end{align*} for any $\Z$-algebra $R$. Then $G_n\times_\Z \Q$ and
$\G_n\times_\Z\Q$ are reductive. We let $\theta$
denote the action on $\G_{ n}$ induced by $(1,c)$ on $G_{n}\times_{\Z} \CO_E$.

If $R$ is an $E$-algebra, then $G_n(R)$ is a subgroup of $R^\times\times
GL_n(F\otimes_\Q R)=R^\times\times GL_n(F\otimes_ER)\times GL_n(F\otimes_{E,c}R)$ and the
projection onto $R^\times\times GL_n(F\otimes_ER)$ defines
an isomorphism
\[ G_n(R) \cong R^\times\times GL_n(F\otimes_ER) .\]
It follows that $G_n\times_\Q E \iso \G_m\times R_{F/E}(GL_n)$.

If $v \in \Unr_{F/\Q}$, then $K_v:=G_n(\Z_v)$ (resp.\ $\K_v:=\G_n(\Z_v)$) is a
hyperspecial maximal compact subgroup of $G_n(\Q_v)$ (resp.\
$\G_n(\Q_v)$). In this case we say that a representation of
$G_n(\Q_v)$ (resp.\ $\G_n(\Q_v)$) is unramified if the space of
$K_v$-invariants (resp.\ $\K_v$-invariants) is
non-zero. Furthermore, we define the unramified Hecke algebras 
$\CH^{\ur}(G_n(\Q_v))$ and $\CH^{\ur}(\G_n(\Q_v))$ with respect to
$K_v$ and $\K_v$ respectively, as in Section 1.1 of
\cite{shin}. (We note that these are $\C$-algebras.) If
$T$ is a set of places of $\Q$ with $\{\infty\}\cup\Ram_{F/\Q}\subset
T$, we let $K^T=\prod_{v\not\in T}K_v\subset G_n(\A^T)$.

We say that a representation $\Pi_v$ of $\G_n(\Q_v)$ is
$\theta$-stable if $\Pi_v\circ \theta \cong \Pi_v$ and we let
$\Irr^{\theta-\st}(\G_n(\Q_v)) \subset \Irr(\G_n(\Q_v))$ be the subset
of $\theta$-stable representations. For $v\in \Unr_{F/\Q}$, we let
$\Irr^{\ur}(G(\Q_v)) \subset \Irr(G(\Q_v))$ (resp.\
$\Irr^{\ur}(\G_n(\Q_v))\subset \Irr(\G_n(\Q_v))$, resp.\
$\Irr^{\ur,\theta-\st}(\G_n(\Q_v))\subset \G_n(\Q_v)$)
denote the subset consisting of 
unramified (resp.\ unramified, resp.\
unramified, $\theta$-stable) representations.

Let $\# : R_{F/\Q}(GL_n) \ra R_{F/\Q}(GL_n)$ denote the map $g\mapsto
\Phi_n{}^t\!g^{-c}\Phi_n^{-1}$. If $v$ is a rational prime, then
\[ \G_n(\Q_v)=G_n(E\otimes_\Q \Q_v) \cong (E\otimes_\Q \Q_v)^\times
\times GL_n(F\otimes_\Q \Q_v). \] 
If $(\lambda,g)\in \G_n(\Q_v)$, then $\theta(\lambda,g)=(\lambda^c,\lambda^cg^\#)$.
Let $\Pi_v \in \Irr(\G_n(\Q_v))$ and write $\Pi_v=\psi_v\otimes\Pi_v^1$ with
respect to the above decomposition of $\G_n(\Q_v)$. Then $\Pi_v$ is
$\theta$-stable if and only if $(\Pi_v^1)^\vee\cong \Pi_v^1\circ c$,
and $\psi_{\Pi_v^1}|_{(E\otimes_\Q\Q_v)^\times}=\psi^c_v/\psi_v$.

We now recall the existence of local base change maps in the
following cases (see Section 4.2 of \cite{shin} for details):

\begin{itemize}
\item  {\it Case 1:} If $v \in \Unr_{F/\Q}$, we have a map \[\BC_v : \Irr^{\ur}(G_n(\Q_v))\ra
\Irr^{\ur,\theta-\st}(\G_n(\Q_v)).\] (Note that the assumption $v \not
\in \Ram_{\Q}(\varpi)$ in Case 1 of Section 4.2 of \cite{shin} plays
no role there.) This is induced by a homomorphism of $\C$-algebras
$\BC_v^* : \CH^{\ur}(\G_n(\Q_v))\ra\CH^{\ur}(G_n(\Q_v))$.

\item  {\it Case 2:} If $v \in \Spl_{F/F^+,\Q}$, we
  have a map 
\[ \BC_v: \Irr(G_n(\Q_v))\ra\Irr^{\theta-\st}(\G_n(\Q_v)).\]
If in addition $v\in \Unr_{F/\Q}$, then this map is compatible with
the map in Case 1.
\end{itemize}

In Case 2, the map $\BC_v$ is described explicitly in
\cite{shin}. We recall the explicit definition here, assuming $v \in
\Spl_{E/\Q}$. Let $y|v$ be a place of $E$, and regard $\Q_v$ as an
$E$-algebra via $\Q_v\iso E_y$. We get an isomorphism
\[ G_n(\Q_v)\isoto \Q_v^\times\times \prod_{x|y}GL_n(F_x) \]
where the product is over all places of $F$ dividing $y$. Let
$\pi_v\in\Irr(G_n(\Q_v))$ and decompose $\pi_v=\pi_{v,0}\otimes\pi_y$
with respect to the above decomposition of $G_n(\Q_v)$. If we
decompose
\[ \G_n(\Q_v)=E_y^\times\times E^\times_{y^c}\times
\prod_{x|y}GL_n(F_x)\times\prod_{x|y}GL_n(F_{x^c}) \]
then $\BC_v(\pi_v)=(\psi_y,\psi_{y^c},\Pi_y,\Pi_{y^c})$, where 
\[
  (\psi_y,\psi_{y^c},\Pi_y,\Pi_{y^c})=(\pi_{p,0},\pi_{p,0}(\psi_{\pi_y}|_{E_y^\times}\circ
  c),\pi_y,\pi_y^\#)
\]
and $\pi_y^{\#}(g):=\pi_y(\Phi_n{}^t\!g^{-c}\Phi_n^{-1})$. (In particular,
$\pi_y^{\#}\cong \pi_y^\vee\circ c$.)

The discussion above can be carried out equally well in the setting
of $\Qlbar$-representations, and we define
$\Irr_l^{\theta-\st}(\G_n(\Q_v))$, $\Irr_l^{\ur}(G(\Q_v))$ etc.\ in
the obvious fashion. We also define a base change map
$\BC_v$ in Case 1 (resp.\ Case 2) by setting $\BC_v(\pi) =
\imath^{-1}\BC_v(\imath \pi)$ for $\pi \in \Irr^{\ur}_l(G_n(\Q_v))$
(resp.\ $\pi \in \Irr(G_n(\Q_v))$).

Let $\xi_\C$ denote an irreducible algebraic representation of $G_n$
over $\C$. There is an isomorphism
$\G_n(\C)=G_n(E\otimes_\Q\C)\cong G_n(\C)\times G_n(\C)$ induced by
the isomorphism $E\otimes_\Q\C\iso \C\times\C$ which sends $e\otimes
z$ to $(\tau(e)z,\tau^c(e)z)$. We associate to $\xi$ a $\theta$-stable irreducible
algebraic representation $\Xi$ of $\G_n$ over $\C$ by setting
$\Xi:=\xi\otimes\xi$. Every such $\Xi$ arises in this way.

We also fix the following data:
\begin{itemize}
\item $V=F^n$ as an $F$-vector space;
\item $\seq{\cdot,\cdot}:V\times V \ra \Q$ is a non-degenerate pairing
  such that $\seq{fv_1,v_2}=\seq{v_1,f^c v_2}$ for all $v_1,v_2 \in V$
  and $f \in F$;
\item $h:\C \ra \End_F(V)\otimes_{\Q}\R$ is an $\R$-algebra embedding
  such that the bilinear pairing $(V\otimes_\Q\R)\times (V\otimes_\Q\R)\ra \R;(v_1,v_2)\mapsto \seq{v_1,h(i)v_2}$
  is symmetric and positive definite. 
\end{itemize}
Under the natural isomorphism $\End_F(V)\otimes_{\Q}\R \cong
\prod_{\sigma\in\Phi^+_\C}M_n(\C)$ we assume that $h$ sends
\[ z \mapsto \left(
    \begin{pmatrix}
      zI_{p_{\sigma}} & 0 \\
      0 & \overline{z}I_{q_{\sigma}}
    \end{pmatrix}
_{\sigma\in\Phi^+_\C}\right)
\]
for some $p_\sigma,q_\sigma \in \Z_{\geq 0}$ with $p_\sigma+q_\sigma=n$.

Define a reductive algebraic group $G/\Q$ by setting
\[ G(R)=\{ (\lambda,g)\in R^\times\times GL_n(F\otimes_\Q R) :
\seq{gv_1,gv_2}=\lambda\seq{v_1,v_2} \mbox{for all } v_1,v_2 \in
V\otimes_\Q R\} \]
for each $\Q$-algebra $R$. Note that $G_n$ is a
quasi-split inner form of $G$. Let $\nu : G\ra \G_m$ denote the
homomorphism which sends $(\lambda,g)$ to $\lambda$.

By Lemma 5.1 of \cite{shin} we can and do assume that
$\seq{\cdot,\cdot}$ and $h$ have been chosen so that 
\begin{itemize}
\item $G_{\Q_v}$ is quasi-split for each rational prime $v$;
\item for each $\sigma\in\Phi_\C^+$, we have
  $(p_\sigma,q_\sigma)=(1,n-1)$ if $\sigma=\tau$ and
  $(p_\sigma,q_\sigma)=(0,n)$ otherwise.
\end{itemize}
As a consequence, we can and do fix an isomorphism
\[ G\times_\Q \A^\infty \cong G_n\times_\Q \A^\infty.\]
Using this isomorphism, we will henceforth identify the groups
$G_n(\Q_v)$ and $G(\Q_v)$ for all primes $v$.
Let $C_G \in \Z_{>0}$ be the integer $|\ker^1(\Q,G)|\cdot \tau(G)$ in
the notation of \cite{shin}.

Let $T$ be a (possibly infinite) set of places of $\Q$ containing
$\infty$ and let $T_{\fin}=T-\{\infty\}$. Let $\Gamma$ be a Galois
group with its Krull topology, or the Weil group of a local field, or
a quotient of such a group. We define an admissible
$\Qlbar[G(\A^T)\times \Gamma]$-module to be an admissible
$\Qlbar[G(\A^T)]$-module $R$ with a commuting continuous action of
$\Gamma$ (the continuity condition here means that for each compact open subgroup
$U\subset G(\A^T)$, the induced map $\Gamma \ra \Aut(R^U)$ is
continuous for the $l$-adic topology on $R^U$). We let $\Groth_l(G(\A^T)\times \Gamma)$ denote the
Grothendieck group of the category of admissible $\Qlbar[G(\A^T)\times
\Gamma]$-modules. If $R$ is an admissible $\Qlbar[G(\A^T)\times
\Gamma]$-module, we let $[R]$ denote its image in
$\Groth_l(G(\A^T)\times\Gamma)$. We let $\Irr_l(G(\A^T)\times\Gamma)$ denote the set
of isomorphism classes of irreducible admissible
$\Qlbar[G(\A^T)\times\Gamma]$-modules. (See Section I.2 of \cite{ht}.)

Now suppose that $T$ is finite, that $p\in T$ and let $J/\Q_p$ be a reductive
group. Let $G'$ be a topological group which is of the form
$G(\A_{T_{\fin}})\times \Gamma$, or
$G(\A_{T_{\fin}-\{p\}})\times\Gamma$, or
$G(\A_{T_{\fin}-\{p\}})\times J(\Q_p)$. Let $[X] \in
\Groth_l(G(\A^T)\times G')$ and write $[X]=\sum_{\pi^T,\rho} n(\pi^T
\otimes \rho)\cdot[\pi^T\otimes \rho]$ where $n(\pi^T\otimes \rho)\in
\Z$ and $\pi^T$ (resp.\ $\rho$) runs through $\Irr_l(G(\A^T))$ (resp.\
$\Irr_l(G')$). For a given $\pi^T \in \Irr_l(G(\A^T))$, we let
\[ [X][\pi^T]= \sum_\rho n(\pi^T\otimes\rho)\cdot[\pi^T\otimes \rho]
\in \Groth_l(G(\A^T)\times G').\]
If $\Ram_{F/\Q}\subset T$ and $\Pi^T\in\Irr(\G_n(\A^T))$ is unramified
at all $v\not\in T$, then we define
\[ [X][\Pi^T] =\sum_{\pi^T} [X][\pi^T] \in \Groth_l(G(\A^T)\times G') 
\]
where the sum is over all $\pi^T\in \Irr_l(G(\A^T))$ with $\pi^T$
unramified at all $v\not\in T$ and
$\BC^T(\imath\pi^T):=\otimes'_{v\not\in T}\BC_v(\imath\pi^T_v)\cong
\Pi^T$.

Finally, suppose $G'$ of the form
$G(\A_{T_{\fin}-\{p\}})\times\Gamma$ or $G(\A_{T_{\fin}})\times
\Gamma$ and let $R$ be an admissible
$\Qlbar[G(\A^T)\times G']$-module. Suppose $\Ram_{F/\Q}\subset T$ and
$\Pi^T\in\Irr(\G_n(\A^T))$ is unramified at all $v\not\in T$. Let
$\CH^{\ur}(G(\A^T))=\otimes'_{v\not\in T} \CH^{\ur}(G(\Q_v))$, a commutative
polynomial algebra over $\C$ in countably many variables. Similarly,
let $\CH^{\ur}(\G_n(\A^T))=\otimes'_{v\not\in T}
\CH^{\ur}(\G_n(\Q_v))$. Then $\Pi^T$ corresponds to a maximal ideal
$\gn$ of $\CH^{\ur}(\G_n(\A^T))$ with residue field $\C$. Note that
the space of $K^T$-invariants $R^{K^T}$ is a module over $\iota^{-1}\CH^{\ur}(G(\A^T))$.
We define
\[ R^{K^T}\{\Pi^T\}:= \bigoplus_{\gm} (R^{K^T})_{\iota^{-1}\gm} \subset R^{K^T} \]
where $\gm$ runs over the maximal ideals of $\CH^{\ur}(G(\A^T))$
with residue field $\C$ and which pull back to $\gn$ under
$\otimes_{v\not\in T}\BC_v^*$. Then
$R^{K^T}\{\Pi^T\}$ is a $G'$-stable direct summand of $R^{K^T}$.

\section{Shimura varieties}
\label{subsec:shimura-varieties}

In this section we recall some results from \cite{shin}. We begin with
some definitions and refer the reader to Section 5 of \cite{shin} for
more details.  Let $U$ be a compact open subgroup of $G(\A^\infty)$ and define a
functor $\gX_U$ from the category of pairs $(S,s)$, where is $S$ is a
connected locally Noetherian $F$-scheme and $s$ is a geometric point
of $S$, to the category of sets by sending a pair $(S,s)$ to the set
of isogeny classes of quadruples $(A,\lambda,i,\bareta)$ where
\begin{itemize}
\item $A/S$ is an abelian scheme of dimension $[F^+:\Q]n$;
\item $\lambda : A\ra A^\vee$ is a polarization;
\item $i:F\into \End(A)\otimes_\Z \Q$ such that $\lambda\circ i(f)=
  i(f^c)^\vee\circ \lambda$;
\item $\bareta$  is a $\pi_1(S, s)$-invariant $U$-orbit of
  isomorphisms of $F\otimes_\Q \A^\infty$-modules $\eta : V\otimes_\Q \A^\infty
  \isoto VA_s$ which take the pairing $\seq{\cdot,\cdot}$ on $V$ to a $(\A^\infty)^\times$-multiple 
   of the $\lambda$-Weil pairing on $VA_s:=H_1(A_s,\A^\infty)$ (see Section 5
  of \cite{kott});
\item for each $f\in F$ there is an equality of polynomials
  $\det_{\CO_S}(f|\Lie A) = \det_E(f|V^1)$ in the sense of Section 5
  of \cite{kott} (here $V^1\subset V\otimes_{\Q}E\subset
  V\otimes_{\Q}\C$ is the $E$-subspace
  where $h(\tau_E(e))$ acts by multiplication by $1\otimes e$ for all
  $e\in E$);
\item two such quadruples $(A,\lambda,i,\bareta)$ and
  $(A',\lambda',i',\bareta')$ are isogenous if there exists an
  isogeny $A\ra A'$ taking $\lambda,i,\bareta$ to $\gamma
  \lambda',i',\bareta'$ for some $\gamma \in \Q^\times$.
\end{itemize}
If $s$ and $s'$ are two geometric points of a connected locally
Noetherian $F$-scheme $S$ then there is a canonical bijection from
$\gX_U(S,s)$ to $\gX_U(S,s')$. We may therefore think of $\gX_U$ as a
functor from connected locally Noetherian $F$-schemes to sets and then
extend it to a functor from all locally Noetherian $F$-schemes to sets
by setting $\gX_U(\coprod_iS_i)=\prod_i \gX_U(S_i)$. When $U$ is
sufficiently small the functor $\gX_U$ is represented by a smooth
projective variety $X_U/F$ of dimension $n-1$. The variety
$X_U$ is denoted $\Sh_U$ in \cite{shin}. Let $\CA_U$ be the universal
abelian variety over $X_U$.

If $U$ and $V$ are sufficiently small compact open subgroups of
$G(\A^\infty)$ and $g\in G(\A^\infty)$ is such that $g^{-1}Vg \subset
U$, then we have a map $g : X_V \ra X_U$ and a quasi-isogeny $g^* :
\CA_V \to g^*\CA_U$ of abelian varieties over $X_V$. In this way we
get a right action of the group $G(\A^\infty)$ on the inverse system of
the $X_U$ which extends to an action by quasi-isogenies on the inverse
system of the $\CA_U$.

Let $l$ be a rational prime and let $\xi$ be an irreducible algebraic representation of $G$ over
$\Qlbar$. Then $\xi$ gives rise to a lisse $l$-adic sheaf $\CL_\xi$ on
each $X_U$. We let
\[ H^k(X,\CL_\xi)= \varinjlim_U H^k(X_U \times_F \barF,\CL_\xi). \]
This is a semisimple admissible representation of $G(\A^\infty)$ with
a commuting continuous action of $G_F$ and therefore decomposes as
\[ H^k(X,\CL_\xi)= \bigoplus_{\pi^\infty} \pi^\infty\otimes
R^k_{\xi,l}(\pi^\infty)\]
where $\pi^\infty$ runs over $\Irr_l(G(\A^\infty))$ and each
$R^k_{\xi,l}(\pi^\infty)$ is a finite dimensional continuous
representation of $G_F$ over $\Qlbar$.

We now recall results of Shin on the existence of Galois
representations in the cohomology of the Shimura varieties in the
following two cases:

\subsection{The stable case}
\label{subsubsec:stable-case}
Assume we are in the following situation:
\begin{itemize}
\item $\Pi^1$ is a RACSDC automorphic representation of $GL_n(\A_F)$;
\item $\psi : \A_E^\times/E^\times\to\C^\times$ is an algebraic Hecke character
  such that $\psi_{\Pi^1}|_{\A_E^\times}=\psi^c/\psi$;
\item $\Pi:=\psi\otimes \Pi^1$ (an automorphic representation of
  $\G_n(\A)\cong GL_1(\A_E)\times GL_n(\A_F)$) is $\Xi$-cohomological for some
  irreducible algebraic representation $\Xi$ of $\G_n/\C$;
\item $\Ram_\Q(\Pi)\subset \Spl_{F/F^+,\Q}$.
\end{itemize}
Then $\Xi$ is $\theta$-stable and so comes from some irreducible algebraic
representation $\xi_\C$ of $G_n$ over $\C$ as in Section \ref{subsec:notat-runn-assumpt}.
Let $\xi=\imath^{-1}\xi_\C$, regarded as a representation of
$G/\Qlbar$. Let $\CR_l(\Pi)$ denote the set of $\pi^\infty \in
\Irr_l(G(\A^\infty))$ such that
\begin{itemize}
\item $R_{\xi,l}^k(\pi^\infty)\neq (0)$ for some $k$;
\item $\pi^\infty$ is unramified at all $v\not \in \Ram_{F/\Q}\cup
  \Ram_\Q(\Pi)$ and 
\[ \BC^\infty(\imath\pi^\infty)=\Pi^\infty. \]
(We note that
$\BC_v(\imath\pi_v)$ is defined for all $v\ndiv \infty$.)
\end{itemize}
 Let $\tR^{k}_{\xi,l}(\Pi)=\bigoplus_{\pi^\infty \in
  \CR_l(\Pi)}R^k_{\xi,l}(\pi^\infty)$.

Now, let $T \supset \{\infty\}$ be a finite set of places of $\Q$ with
$\Ram_{F/\Q}\cup\Ram_{\Q}(\Pi)\subset T_{\fin}\subset
\Spl_{F/F^+,\Q}$. Note that
\[ H^k(X,\cL_\xi)^{K^T}\{\Pi^T\}\cong \bigoplus_{\pi^{\infty}}\pi_{T_{\fin}}\otimes
R_{\xi,l}^k(\pi^\infty) \subset H^k(X,\cL_\xi)\]
where the sum is over all $\pi^\infty=\pi^T\otimes \pi_{T_{\fin}}$
where $\pi^T$ is unramified and $\BC(\imath\pi^T)\cong \Pi^T$. We then
define an admissible $\Qlbar[\G_n(\A_{T_{\fin}})\times G_F]$-module
\[ \BC_{T_{\fin}}(H^k(X,\cL_\xi)^{K^T}\{\Pi^T\}):=\bigoplus_{\pi^{\infty}}\BC_{T_{\fin}}(\pi_{T_{\fin}})\otimes
R_{\xi,l}^k(\pi^\infty), \]
where $\pi^\infty$ runs over the same set.

\begin{thm}
  \label{thm:cohomology-stable-case}
  \begin{enumerate}
  \item If $\pi^\infty\in \CR_l(\Pi)$ then $R^k_{\xi,l}(\pi^\infty)\neq
    (0)$ if and only if $k=n-1$.
  \item We have
\[ \BC_{T_{\fin}}(H^{n-1}(X,\cL_\xi)^{K^T}\{\Pi^T\})=(\imath^{-1}\Pi_{T_{\fin}})\otimes
\tR_{\xi,l}^{n-1}(\Pi).\]
 \item We have  
\[ \tR_{\xi,l}^{n-1}(\Pi)^{\semis} \cong R_{l,\imath}(\Pi^1)^{C_G}\otimes
R_{l,\imath}(\psi)|_{G_F}. \]
  \end{enumerate}
\end{thm}

\begin{proof}
  The first part follows from Corollary 6.5 of \cite{shin}. The second
  part follows from the proof of Corollary 6.4 of {\it op.\ cit.}. The third
  part follows from the proof of Corollary 6.8 of {\it op.\ cit.}
  (Note that the character $\rec_{l,\imath_l}(\psi)$ which appears in the
  proof of this corollary is equal to
  $R_{l,\imath}(\psi^{-1})$.)
\end{proof}

\subsection{The endoscopic case}
\label{subsubsec:endoscopic-case}

We now assume we are in the following situation:
\begin{itemize}
\item $m_1,m_2$ are positive integers with $m_1>m_2$ and $m_1+m_2=n$;
\item for $i=1,2$, $\Pi_i$ is a RACSDC automorphic representation of
  $GL_{m_i}(\A_F)$ with $\Ram_\Q(\Pi_i)\subset \Spl_{F/F^+,\Q}$;
\item $\varpi: \A_E^\times/E^\times \ra \C^\times$ is a Hecke character such
  that
  \begin{itemize}
  \item $\Ram_{\Q}(\varpi)\subset \Spl_{F/F^+,\Q}$;
  \item $\varpi|_{\A^\times_\Q/\Q^\times}$ is the quadratic character
    corresponding to the quadratic extension $E/\Q$ by class field
    theory.
  \end{itemize}
\item $\psi: \A_E^\times/E^\times\ra \C^\times$ is an algebraic Hecke character
  such that
  \begin{itemize}
  \item $(\psi_{\Pi_1}\psi_{\Pi_2})|_{\A_E^\times}=(\psi\otimes
    \varpi^{N(m_1,m_2)})^c/(\psi\otimes\varpi^{N(m_1,m_2)})$ where $N(m_1,m_2)=[F^+:\Q](m_1\epsilon(n-m_1)+m_2\epsilon(n-m_1))/2\in\Z$;
  \item $\Ram_{\Q}(\psi) \subset \Spl_{F/F^+,\Q}$.
  \end{itemize}
\item for $i=1,2$, $\Pi_{M,i}:=\Pi_i\otimes (\varpi\circ \norm_{F/E}\circ
  \det)^{\epsilon(n-m_i)}$;
\item $\Pi:=\psi\otimes \nind_{GL_{m_1}\times GL_{m_2}}^{GL_n}(\Pi_{M,1}\otimes
  \Pi_{M,2})$ (an automorphic representation of $\G_n(\A)$) is
  $\Xi$-cohomological for some irreducible algebraic representation
  $\Xi$ of $\G_n/\C$. (We note that the normalized induction is
  irreducible as $\Pi_{M,1}$ and $\Pi_{M,2}$ are unitary.)
\end{itemize}
As above, we let $\xi$ be the irreducible algebraic
representation of $G$ over $\Qlbar$ such that $\Xi$ is associated to
$\imath\xi$. Let $\CR_l(\Pi)$ denote the set of $\pi^\infty \in
\Irr_l(G(\A^\infty))$ such that
\begin{itemize}
\item $R_{\xi,l}^k(\pi^\infty)\neq (0)$ for some $k$;
\item $\pi^\infty$ is unramified at all $v\not \in \Ram_{F/\Q}\cup
  \Ram_\Q(\Pi)\cup\Ram_{\Q}(\varpi)$ and 
\[ \BC^\infty(\imath\pi^\infty)=\Pi^\infty. \]
\end{itemize}
Let $\tR^{k}_{\xi,l}(\Pi)=\bigoplus_{\pi^\infty \in \CR_l(\Pi)}R^k_{\xi,l}(\pi^\infty)$.

Let $T \supset \{\infty\}$ be a finite set of places of $\Q$ with
$\Ram_{F/\Q}\cup\Ram_{\Q}(\Pi)\cup\Ram_{\Q}(\varpi)\subset T_{\fin}\subset
\Spl_{F/F^+,\Q}$. We define
\[ \BC_{T_{\fin}}( H^k(X,\cL_\xi)^{K^T}\{\Pi^T\}) \]
exactly as in the previous subsection. 

\begin{thm}
  \label{thm:cohomology-endo-case}
  \begin{enumerate}
  \item If $\pi^\infty\in \CR_l(\Pi)$, then
    $R^k_{\xi,l}(\pi^\infty)\neq (0)$ if and only if $k=n-1$.
  \item We have
\[ \BC_{T_{\fin}}( H^k(X,\cL_\xi)^{K^T}\{\Pi^T\})=(\imath^{-1}\Pi_{T_{\fin}})\otimes
\tR_{\xi,l}^{n-1}(\Pi).\]
  \item There exists an integer $e_2(\Pi,G) \in \{\pm 1\}$ depending on $\Pi$ and $G$ such that
    \begin{enumerate}
    \item\label{a} If $e_2(\Pi,G)=1$ then
      \[ \tR^{n-1}_{\xi,l}(\Pi)^{\semis} \cong
      R_{l,\imath}(\Pi_1)^{C_G}\otimes
      R_{l,\imath}(\psi\varpi^{\epsilon(n-m_1)}|\cdot|^{(n-m_1)/2})|_{G_F}. \]
    \item\label{b} If $e_2(\Pi,G)=-1$ then
\[ \tR^{n-1}_{\xi,l}(\Pi)^{\semis}  \cong R_{l,\imath}(\Pi_2)^{C_G}\otimes
R_{l,\imath}(\psi\varpi^{\epsilon(n-m_2)}|\cdot|^{(n-m_2)/2})|_{G_F}. \]
    \end{enumerate}
  \end{enumerate}
\end{thm}

\begin{proof}
  The first and second parts follow respectively from Corollary 6.5
  and the proof of Corollary 6.4 of \cite{shin}. The
  third part follows from the proofs of Corollaries 6.8 and 6.10 of {\it
    op.\ cit.} (Alternative \ref{a} corresponds to the case when
  $e_1=e_2$ in the notation of {\it op. cit.}, while alternative \ref{b}
  corresponds to the case when $e_1=-e_2$. Note however that by
  Corollary 6.5 of {\it op. cit.}, $e_1=(-1)^{n-1}=1$. We therefore
  take $e_2(\Pi,G)=e_2$.)
\end{proof}

\section{Integral models}
\label{subsec:integral-models}
We now proceed to introduce integral models for the varieties
$X_U$ and to deduce various results on these models, following the arguments of Section 3 of \cite{ty}.
Recall that we have fixed an isomorphism
$G\times_{\Q}\A^\infty \cong G_n\times_{\Q}\A^\infty$. Since $p \in
\Spl_{E/\Q}$, we have an isomorphism
\[ G(\Q_p) \cong \Q_p^\times\times\prod_{i=1}^r GL_n(F_{w_i}) \]
and we decompose $G(\A^\infty)$ as
\[ G(\A^\infty)=G(\A^{\infty,p})\times\Q_p^\times\times\prod_{i=1}^r GL_n(F_{w_i}).\] 
If $m=(m_2,\ldots,m_r)\in \Z^{r-1}_{\geq 0}$, set
\[ U_p^w(m) = \prod_{i=2}^r \ker(GL_n(\CO_{F,w_i})\ra GL_n(\cO_{F,w_i}/w_i^{m_i}))
\subset \prod_{i=2}^r GL_n(F_{w_i}).\]
We consider the following compact open subgroups of $G(\Q_p)$:
\begin{align*}
  \Ma(m)&=\Z_p^\times \times GL_n(\CO_{F,w}) \times U_p^w(m) \\
  \Iw(m)&=\Z_p^\times \times \Iw_{n,w} \times U_p^w(m) .
\end{align*}
Fix an $m$ as above. If $U^p \subset G(\A^{\infty,p})$ is a compact open subgroup, we let
$U_0 = U^p\times \Ma(m)$ and $U=U^p \times \Iw(m)$. For
$i=1,\ldots,r$, let $\Lambda_i \subset V\otimes_F F_{w_i}$ be a
$GL_n(\cO_{F_{w_i}})$-stable lattice.

For each sufficiently small $U^p$ as above, an integral model of $X_{U_0}$ over
$\CO_{F,w}$ is constructed in Section 5.2 of \cite{shin} (note that
$X_{U_0}$ is denoted $\Sh_{U^p(\vec m)}$ in \cite{shin} with $\vec m=(0,m_2,\ldots,m_r)$). We denote
this integral model also by $X_{U_0}$. It represents a functor
$\gX_{U_0}$ from the category of locally Noetherian
$\CO_{F,w}$-schemes to sets which, as in the characteristic 0 case, is
initially defined on the category of pairs $(S,s)$ where $S$ is a connected locally Noetherian
$\CO_{F,w}$-scheme and $s$ is a geometric point of $S$. It sends a
pair $(S,s)$ to the set of equivalence classes of $(r+3)$-tuples $(A,\lambda,i,\bareta^p,\{\alpha_i\}_{i=2}^r)$
where
\begin{itemize}
\item $A/S$ is an abelian scheme of dimension $[F^+:\Q]n$;
\item $\lambda : A \ra A^\vee$ is a prime-to-$p$ polarization;
\item $i:\CO_F \into \End(A)\otimes_\Z\Z_{(p)}$ such that
  $\lambda\circ i(f)=i(f^c)^\vee \circ \lambda$;
\item $\bareta^p$  is a $\pi_1(S, s)$-invariant $U^p$-orbit of
  isomorphisms of $F\otimes_\Q \A^{\infty,p}$-modules $\eta^p : V\otimes_\Q \A^{\infty,p}
  \isoto V^p A_s$ which take the pairing $\seq{\cdot,\cdot}$ on $V$ to a $(\A^{\infty,p})^\times$-multiple 
   of the $\lambda$-Weil pairing on $V^p A_s$;
\item for each $f\in \CO_F$ there is an equality of polynomials
  $\det_{\CO_S}(f|\Lie A) = \det_E(f|V^1)$ in the sense of Section 5
  of \cite{kott};
\item for $2\leq i\leq r$, $\alpha_i:
  (w_i^{-m_i}\Lambda_i/\Lambda_i)\isoto A[w_i^{m_i}]$ is an
  isomorphism of $S$-schemes with $\CO_{F,w_i}$-actions;
\end{itemize}
and
\begin{itemize}
\item two such tuples $(A,\lambda,i,\bareta^p,\{\alpha_i\}_{i=2}^r)$
  and $(A',\lambda',i',(\bareta^p)',\{\alpha'_i\}_{i=2}^r)$ are
  equivalent if there is a prime-to-$p$ isogeny $A\ra A'$ taking
  $\lambda$, $i$, $\bareta^p$ and $\alpha_i$ to
  $\gamma\lambda'$, $i'$, $(\bareta^p)'$ and $\alpha_i'$ for some $\gamma\in\Z_{(p)}^\times$. 
\end{itemize}
The scheme $X_{U_0}$ is smooth and projective over $\CO_{F,w}$. As
$U^p$ varies, the inverse system of the $X_{U_0}$'s has an action of
$G(\A^{\infty,p})$.

Given a tuple $(A,\lambda,i,\bareta^p,\{\alpha_i\}_{i=2}^r)$ over $S$ as above,
we let $\CG_A=A[w^\infty]$, a Barsotti-Tate $\CO_{F,w}$-module over $S$. If $p$
is locally nilpotent on $S$, then $\CG_A$ has dimension 1 and is
compatible (which means that the two actions of $\CO_{F,w}$ on $\Lie
\CG_A$ (coming from the structural morphism $S\ra \Spec \cO_{F,w}$
and from $i:\cO_F \ra \End(A)\otimes_\Z\Z_{(p)}$) coincide). We let
$\CA_{U_0}$ denote the universal abelian scheme over $X_{U_0}$, and we
let $\CG=\CG_{\CA_{U_0}}$.

Let $\barX_{U_0}$ denote the special fibre $X_{U_0}\times_{\CO_{F,w}}k(w)$
of $X_{U_0}$, and for $0\leq h \leq n-1$, let $\barX_{U_0}^{[h]}$ denote the reduced closed
subscheme of $\barX_{U_0}$ whose closed geometric points $s$
are those for which the maximal \'etale quotient of $\CG_s$ has
$\CO_{F,w}$-height at most $h$. Let
\[ \barX_{U_0}^{(h)}=\barX_{U_0}^{[h]}-\barX_{U_0}^{[h-1]} \]
(where we set $\barX^{[-1]}_{U_0}=\emptyset$). Then
$\barX_{U_0}^{(0)}$ is non-empty. [We exhibit an $\barFF_p$ point of $\barX_{U_0}^{(0)}$. Consider the $p$-adic type $(F,\eta)$ over $F$ where $\eta_w=1/(n [k(w):\F_p])$ and $\eta_{w_i}=0$ for $i>1$. It corresponds to an isogeny class of abelian varieties with $F$-action over $\barFF_p$. Let $(A,i)/\barFF_p$ be an element of this isogeny class. Then 
\begin{itemize}
\item $A$ has dimension  $[F^+:\Q]n$; 
\item $C_0=\End^0_F(A)$ is the division algebra with centre $F$ which is split outside $ww^c$ and has Hasse invariant $1/n$ at $w$; 
\item and the $p$-divisible group $A[w^\infty]$ has pure slope $1/(n[F_w:\Q_p])$, while $A[w_i^\infty]$ is \'{e}tale for $i>1$. 
\end{itemize}
(See section 5.2 of \cite{ht}.) Just as in the proof of Lemma V.4.1 of \cite{ht} one shows that there is a
polarization $\lambda_0:A \ra A^\vee$ and an $F$-vector space $W_0$ of dimension $n$ together with a non-degenerate alternating form $\langle \,\,\, , \,\,\, \rangle_0: W_0 \times W_0 \ra \Q$ such that 
\begin{itemize}
\item $\lambda_0 \circ i(a) = i(c a)^\vee \circ \lambda_0$ for all $a \in F$;
\item $\langle ax,y \rangle_0=\langle x,(ca)y \rangle_0$ for all $a \in F$ and $x,y \in W_0$;
\item $V^pA \cong W_0 \otimes_\Q \A^{\infty,p}$ as $\A_F^{\infty,p}$-modules with alternating pairings defined up to $(\A^{\infty,p})^\times$-multiples (the pairing on $V^pA$ being the $\lambda_0$-Weil pairing);
\item $W_0 \otimes_\Q \R \cong V \otimes_\Q \R$ as $F \otimes_\Q \R$-modules with alternating pairings up to $\R^\times$-multiples.
\end{itemize}
(In fact $W_0$ will be the Betti cohomology of a certain lift of $(A,i)$ to characteristic $0$.)
Let $G_0$ denote the denote the algebraic group of $F$-linear automorphisms of $W_0$ that preserve
$\langle \,\,\, ,\,\,\, \rangle_0$ up to scalar multiples, and let $\phi_0 \in H^1(\Q, G_0)$ represent the difference between $(W_0, \langle\,\,\, ,\,\,\,\rangle_0)$ and $(V,\langle \,\,\, ,\,\,\, \rangle)$. So in fact 
\[ \phi_0 \in \ker (H^1(\Q,G_0) \lra H^1(\R,G_0)). \]
Let $\ddag_0$ denote the $\lambda_0$ Rosati involution on $C_0$ and define an algebraic group $H_0^\AV/\Q$ by
\[ H_0^\AV(R)=\{ g \in (C_0 \otimes_\Q R)^\times : g g^{\ddag_0} \in R^\times \}. \]
There is a natural isomorphism
\[ H_0^\AV \times_\Q \A^{\infty,p} \liso G_0 \times_\Q \A^{\infty,p} \]
coming from the isomorphism  $V^pA \cong W_0 \otimes_\Q \A^{\infty,p}$. As in Lemma V.3.1 of \cite{ht} the polarizations of $A$ which induce complex conjugation on $i(F)$ are parametrized by 
\[ \ker (H^1(\Q,H_0^\AV) \lra H^1(\R,H_0^\AV)). \]
Let $A(G_0)$ and $A(H_0^\AV)$ be the groups defined in Section 2.1 of \cite{kotstfest}, so that there are
sequences
\[ H^1(\Q,G_0) \lra H^1(\Q,G_0(\barAA)) \lra A(G_0) \]
and
\[ H^1(\Q,H^\AV_0) \lra H^1(\Q,H^\AV_0(\barAA)) \lra A(H^\AV_0), \]
which are exact in the middle.
Note that as all primes of $F^+$ above $p$ split in $F$ we have $H^1(\Q_p,G_0)=(0)$ and $H^1(\Q_p,H_0^\AV)=(0)$. Thus the image $\phi_0^{\infty,p}$ of $\phi_0$ in $H^1(\Q,G_0(\barAA^{\infty,p}))$ maps to $0$ in $A(G_0)$. By Lemma 2.8 of \cite{kott} 
\[ \begin{array}{ccc} H^1(\Q,H_0^\AV(\barAA^{\infty,p})) & \liso & H^1(\Q,G_0(\barAA^{\infty,p})) \\ \da & & \da \\ A(H_0^\AV) & = & A(G_0) \end{array} \]
commutes. Thus thinking of $\phi_0^{\infty,p} \in H^1(\Q,H_0^\AV(\barAA^{\infty,p}))$ we see that it can
be lifted to 
\[ \phi_0^\AV \in \ker(H^1(\Q,H_0^\AV) \lra H^1(\R,H_0^\AV)). \]
Let $\lambda$ denote the corresponding polarization of $(A,i)$. There is an isomorphism
\[ \eta^p: V \otimes_\Q \A^{\infty,p} \liso V^pA \]
of $\A_F^{\infty,p}$-modules with alternating pairings up to $(\A^{\infty,p})^\times$-multiples. Moroever for $i=2,...,r$ the $p$-divisible group $A[w_i^\infty]$ is \'{e}tale and so there are isomorphisms $(w_i^{-m_i}\Lambda_i /\Lambda_i \iso A[w_i^{m_i}]$. Thus 
\[ (A,\lambda,i,\bareta^p,\{ \alpha_i\}_{i=2}^r) \in \barX_{U_0}^{(0)}(\barFF_p), \]
as desired.]

Just as in Section III.4 of \cite{ht}, one deduces that each $\barX_{U_0}^{(h)}$  is non-empty and smooth of pure
dimension $h$.
Over $\barX_{U_0}^{(h)}$ there is a short exact sequence
\[ 0 \ra \CG^0\ra\CG\ra\CG^\et\ra 0\]
where $\CG^0$ is a formal Barsotti-Tate $\CO_{F,w}$-module and
$\CG^{\et}$ is an \'etale Barsotti-Tate $\CO_{F,w}$-module of
$\CO_{F,w}$-height $h$.

We define an integral model for $X_U$
over $\CO_{F,w}$ (for sufficiently small $U^p$) as on page 480 of \cite{ty}. It represents a functor
$\gX_{U}$ from the category of locally Noetherian
$\CO_{F,w}$-schemes to sets which, as above, is
initially defined on the category of connected locally Noetherian
$\CO_{F,w}$-schemes with a geometric point. It sends a
pair $(S,s)$ to the set of equivalence classes of $(r+4)$-tuples $(A,\lambda,i,\bareta^p,\CC,\{\alpha_i\}_{i=2}^r)$
where $(A,\lambda,i,\bareta^p,\{\alpha_i\}_{i=2}^r)$ is as in the
definition of $\gX_{U_0}(S,s)$ and $\CC$ is a chain of isogenies 
\[ \CC : \CG_A=\CG_0\ra\CG_1\ra\cdots\ra\CG_n=\CG_A/\CG_A[w] \]
of compatible Barsotti-Tate $\CO_{F,w}$-modules, each of degree
$\#k(w)$ and with composition equal to the canonical map
$\CG_A\ra\CG_A/\CG_A[w]$. By Lemma 3.2 of \cite{ty}, which holds
equally well in our situation, the functor $\gX_U$ is representable by
scheme $X_U$ which is finite over $X_{U_0}$.

Let $U^p$ be sufficiently small and let
$\barX_U=X_U\times_{\CO_{F,w}}k(w)$ denote the special fibre of
$X_U$. By parts (1) and (2) of Proposition 3.4 of \cite{ty} (whose
proof applies in our situation), $X_U$ has
pure dimension $n$, it has semistable reduction over $\CO_{F,w}$, it
is regular and the natural map $X_U \ra X_{U_0}$ is finite and
flat. We let $\CA_{U}$ denote the universal abelian variety over $X_U$.

We say that an isogeny $\CG \ra \CG'$ of one-dimensional
compatible Barsotti-Tate $\CO_{F,w}$-modules of degree $\#k(w)$ over a
scheme $S$ of characteristic $p$ has {\it connected kernel} if it
induces the zero map $\Lie\CG\ra\Lie\CG'$. Let $Y_{U,i}$ denote the
closed subscheme of $\barX_{U}$ over which $\CG_{i-1}\ra\CG_i$ has
connected kernel. By part (3) of Proposition 3.4 of \cite{ty}, each
$Y_{U,i}$ is smooth over $\Spec k(w)$ of pure dimension $n-1$,
$\barX_U = \cup_{i=1}^n Y_{U,i}$ and for $i\neq j$ the schemes
$Y_{U,i}$ and $Y_{U,j}$ have no common connected component. It follows
that $X_U$ has strictly semistable reduction.

For each $S\subset \{1,\ldots,n\}$, we let
\[ Y_{U,S}=\bigcap_{i\in S}Y_{U,i}\mbox{\ and\ } Y^0_{U,S} =
Y_{U,S}-\bigcup_{T\supsetneq S} Y_{U,T}. \] 
Since $X_U$ has strictly semistable reduction, each $Y_{U,S}$ is
smooth over $k(w)$ of pure dimension $n-\#S$ and the $Y^0_{U,S}$ are
disjoint for different $S$.

The inverse systems $X_U$ and $X_{U_0}$, for varying $U^p$, have compatible
actions of $G(\A^{\infty,p})$. For each $S\subset \{1,\ldots,n\}$, the
inverse systems $Y_{U,S}$ and $Y^0_{U,S}$ are stable
under this action. As in the characteristic zero case, the actions of
$G(\A^{\infty,p})$ extend to actions on the inverse systems of the
universal abelian varieties $\CA_U$ and $\CA_{U_0}$. Here the action
is by prime-to-$p$ quasi-isogenies. 

Let $\xi$ be an irreducible algebraic representation of
$G$ over $\Qlbar$. If $l\neq p$, then the sheaf $\CL_\xi$ extends to a lisse sheaf on the integral
models $X_U$ and $X_{U_0}$. There exist non-negative integers $m_\xi$ and
$t_\xi$ and an idempotent $\varepsilon_\xi\in \Qlbar[S_{m_\xi} \ltimes
F^{m_\xi}]$ (where
$S_{m_\xi}$ is the symmetric group on $m_\xi$-letters) such that
\[ \xi \cong \nu^{t_\xi} \otimes \varepsilon_\xi(V^\vee
\otimes_{\Q}\Qlbar)^{\otimes m_\xi}.\]
This follows from the discussion on pages 97 and 98 of
\cite{ht} (applied in our setting). Let $N\geq 2$ be
prime to $p$ and let
\[ \varepsilon(m_\xi,N) = \prod_{x=1}^{m_{\xi}}\prod_{y\neq 1}
\frac{[N]_x-N^y}{N-N^y} \in\Q[(N^{\Z_{\geq 0}})^{m_\xi}]\]
where $[N]_x$ is the element of $(N^{\Z_{\geq 0}})^{m_\xi}$ with $N$ in the
$x$-th entry and 1 in the other entries and where $y$ runs from 0 to
$2[F^+:\Q]n$ but excludes 1. Thinking of $(N^{\Z_{\geq
    0}})^{m_\xi}\subset F^{m_\xi}$, we set
\[
a_\xi=a_{\xi,N}=\varepsilon_\xi\varepsilon(\xi,N)^{2n-1}\in\Qlbar[S_{m_\xi}\ltimes
F^{m_\xi}].\]
Let $\proj$ denote the map $\CA_U^{m_\xi}\ra X_U$ and for
$a=1,\ldots,m_\xi$, let $\proj_a:\CA_U\ra X_U$ denote the composition of the $a$-th
inclusion $\CA_U\into\CA_U^{m_\xi}$ with $\proj$.
Then we have the following (see page 477 of \cite{ty}):
\begin{itemize}
\item $\varepsilon(\xi,N)R^j\proj_*\Qlbar=\begin{cases} (0) & (j\neq
      m_\xi) \\
\bigotimes_{a=1}^{m_\xi}R^1\proj_{a,*}\Qlbar & (j=m_\xi)\end{cases}$;
\item $\varepsilon_\xi\varepsilon(\xi,N)R^{m_\xi}\proj_*\Qlbar=\CL_\xi$;
\item $a_\xi$ acts as an idempotent on each
$H^j(\CA_U^{m_\xi}\times_{F_w}\barF_w,\Qlbar(t_\xi))$  and
moreover
\[a_\xi H^j(\CA_U^{m_\xi}\times_{\CO_{F,w}}\barF_w,\Qlbar(t_\xi)) =  
 \begin{cases}
   (0) & (j< m_\xi) \\
   H^{j-m_\xi}(X_U \times_{F_w}\barF_w,\CL_\xi) & (j\geq m_\xi).
 \end{cases}
 \]
 \end{itemize}

 Let $\CA_{U,S}^{m_{\xi}}=\CA_U^{m_\xi} \times_{X_U}Y_{U,S}$. As $U^p$
 varies, the inverse system of the $\CA_{U,S}^{m_\xi}$ inherits an action of
$G(\A^{p,\infty})$ by prime-to-$p$ quasi-isogenies. We now make the
following definitions. 
\begin{itemize}
\item Define admissible $\Qlbar[G(\A^{\infty,p})\times G_F]$-modules:
  \begin{align*}
    H^j(X_{\Iw(m)},\CL_\xi)&:=\varinjlim_{U^p}H^j(X_U\times_F\barF,\CL_\xi)=H^j(X,\CL_\xi)^{\Iw(m)} \\
    H^j(\CA_{\Iw(m)}^{m_\xi},\Qlbar)&:=\varinjlim_{U^p}H^j(\CA_U^{m_\xi}\times_F\barF,\Qlbar).
  \end{align*}
\item If $l\neq p$, we define admissible $\Qlbar[G(\A^{p,\infty})\times \Frob_w^\Z]$-modules:
  \begin{align*}
    H^j(Y_{\Iw(m),S},\CL_\xi)&:=\varinjlim_{U^p}H^j(Y_{U,S}\times_{k(w)}\overline{k(w)},\CL_\xi) \\
    H^j_c(Y^0_{\Iw(m),S},\CL_\xi)&:=\varinjlim_{U^p}H^j_c(Y^0_{U,S}\times_{k(w)}\overline{k(w)},\CL_\xi) \\
    H^j(\CA_{\Iw(m),S},\CL_\xi)&:=\varinjlim_{U^p}H^j(\CA_{U,S}\times_{k(w)}\overline{k(w)},\Qlbar) .
  \end{align*}
\item If $l=p$ and $\sigma : W_0 \into \Qlbar$ over $\Z_p=\Z_l$, where
  $W_0$ is the Witt  ring of $k(w)$, then we define the admissible $\Qlbar[G(\A^{\infty,p})\times\Frob_w^\Z]$-module
\[ H^j(\CA_{\Iw(m),S}/W_0)\otimes_{W_0,\sigma}\Qlbar :=
\varinjlim_{U^p} H^j(\CA_{U,S}/W_0)\otimes_{W_0,\sigma}\Qlbar. \]
(Here $H^j(\CA_{U,S}/W_0)$ denotes crystalline cohomology and
$\Frob_w$ acts by the $[k(w):\F_p]$-power of the crystalline Frobenius.)
\end{itemize}
We note that if $l\neq p$, then $a_\xi$ acts as an idempotent on
$H^j(\CA_{\Iw(m),S},\Qlbar)$ and
\[ a_\xi H^j(\CA_{\Iw(m),S},\Qlbar) =
\begin{cases}
  (0) & (j<m_\xi)\\
H^{j-m_\xi}(Y_{\Iw(m),S},\CL_\xi) & (j\geq m_\xi).
\end{cases}\]
If $l=p$, then $a_\xi$ acts as an idempotent on
$H^j(\CA_{\Iw(m),S}/W_0)\otimes_{W_0,\sigma}\Qlbar$. We also note that $I_{F_w}$ acts trivially on
$\WD(H^j(X_{\Iw(m)},\CL_\xi)|_{G_{F_w}})$ and thus the latter can
be regarded as a $\Frob_w^\Z$-module.

\begin{prop}
\label{prop:RZ-spec-seq}
  Let $T$ be a finite set of places of $\Q$ containing
  $\{p,\infty\}\cup \Ram_{F/\Q}$
  and let $\Pi^T\in\Irr(\G_n(\A^T))$ be unramified at all $v\not\in
  T$.  If $l=p$, let $\sigma:W_0 \into \Qlbar$ over $\Z_l=\Z_p$. Then
  there is a spectral sequence in the category of admissible
  $\Qlbar[G(\A_{T_{\fin}-\{p\}})\times \Frob_w^{\Z}]$-modules
\[ E^{i,j}_1(\Iw(m),\xi)^{K^T}\{\Pi^T\} \implies
\WD(H^{i+j}(X_{\Iw(m)},\CL_\xi)|_{G_{F_w}})^{K^T}\{\Pi^T\} \]
where $E^{i,j}_1(\Iw(m),\xi)=\bigoplus_{s\geq
  \max(0,-i)}\bigoplus_{\#S=i+2s+1}H^j_{S,s}$, and
\[ H^j_{S,s}=
\begin{cases}
  a_\xi H^{j+m_\xi-2s}(\CA_{\Iw(m),S}^{m_\xi},\Qlbar(t_\xi-s))=
  H^{j-2s}(Y_{\Iw(m),S},\CL_\xi(-s))& (l\neq p) \\
  a_\xi
  H^{j+m_\xi-2s}(\CA_{\Iw(m),S}^{m_\xi}/W_0)\otimes_{W_0,\sigma}\Qlbar(t_\xi-s)
  & (l=p).
\end{cases}\]
Moreover, the monodromy operator $N$ on
$\WD(H^{i+j}(X_{\Iw(m)},\CL_\xi)|_{G_{F_w}})^{K^T}\{\Pi^T\}$
is induced by the identity map
\begin{gather*} N : \bigoplus_{\#S=i+2s+1}a_\xi
H^{j+m_\xi-2s}(\CA_{\Iw(m),S}^{m_\xi},\Qlbar(t_\xi-s)) \\
\isoto \bigoplus_{\#S=(i+2)+2(s-1)+1}a_\xi
H^{(j-2)+m_\xi-2(s-1)}(\CA_{\Iw(m),S}^{m_\xi},\Qlbar(t_\xi-(s-1)))\end{gather*}
in the case when $l\neq p$
(resp.
\begin{gather*}  N : \bigoplus_{\#S=i+2s+1}a_\xi
 H^{j+m_\xi-2s}(\CA_{\Iw(m),S}^{m_\xi}/W_0)\otimes_{W_0,\sigma}\Qlbar(t_\xi-s)) \\
\isoto \bigoplus_{\#S=(i+2)+2(s-1)+1}a_\xi
  H^{(j-2)+m_\xi-2(s-1)}(\CA_{\Iw(m),S}^{m_\xi}/W_0)\otimes_{W_0,\sigma}\Qlbar(t_\xi-(s-1)) \end{gather*}
in the case when $l=p$).
\end{prop}

\begin{proof}
The proof of Proposition 3.5 of \cite{ty} shows that we have a
spectral sequence $E^{i,j}_1(\Iw(m),\xi) \implies
\WD(H^{i+j}(X_{\Iw(m)},\CL_\xi)|_{G_{F_w}})$ and that the monodromy
operator $N$ on $\WD(H^{i+j}(X_{\Iw(m)},\CL_\xi)|_{G_{F_w}})$ is
induced by the maps above. The result now follows from the fact that
$R\mapsto R^{K^T}\{\Pi^T\}$ is an exact functor from the category of
admissible $\Qlbar[G(\A^{\infty,p})\times \Frob_w^{\Z}]$-modules to
the category of admissible $\Qlbar[G(\A_{T_{\fin}-\{p\}})\times \Frob_w^\Z]$-modules.
\end{proof}

\section{Relating the cohomology of $Y_{U,S}$ to the cohomology of
  Igusa varieties}
\label{sec:rel-cohom-Y}

Let $U^p \subset G(\A^{\infty,p})$ be sufficiently small and let $m\in
\Z^{r-1}_{\geq 0}$. Following Section 4 of \cite{ty}, we can relate the cohomology of
the open strata $Y^0_{U,S}$ to the cohomology of Igusa varieties of
the first kind. For $h=0,\ldots,n-1$ and $m_1 \in \Z_{\geq 0}$, let $I^{(h)}_{U^p,(m_1,m)}/\barX^{(h)}_{U_0}$
denote the Igusa variety of the first kind defined as on page 121 of \cite{ht}.
It is the moduli space of isomorphisms
\[
\alpha_1^{\et}:(w^{-m_1}\CO_{F,w}/\CO_{F,w})^h_{\barX^{(h)}_{U_0}}\isoto
\CG^{\et}[w^{m_1}].\] Let $I^{(h)}_U/\barX^{(h)}_{U_0}$ be the
Iwahori-Igusa variety of the first kind defined as on page 487 of
\cite{ty}. It is the moduli space of chains of isogenies
\[ \CG^{\et}=\CG_0\ra\CG_1\ra\cdots\ra\CG_h=\CG^\et/\CG^{\et}[w]\] of
\'etale Barsotti-Tate $\CO_{F,w}$-modules, each of degree $\#k(w)$ and
with composition equal to the natural map
$\CG^{\et}\ra\CG^{\et}/\CG^{\et}[w]$. Then $I^{(h)}_{U^p,(m_1,m)}$ and
$I^{(h)}_U$ are finite \'etale over $\barX^{(h)}_{U_0}$ and the
natural map $I^{(h)}_{U^p,(1,m)}\ra I^{(h)}_U$ is finite \'etale and
Galois with Galois group $B_h(k(w))$. The inverse systems
$I^{(h)}_{U^p,(m_1,m)}$ and
$I^{(h)}_U$, for varying $U^p$, inherit an action of
$G(\A^{\infty,p})$. Let $\xi$ be an irreducible algebraic
representation of $G$ over $\Qlbar$. If $l\neq p$, then $\xi$ gives
rise to a lisse sheaf $\CL_\xi$ on $I^{(h)}_{U^p,(m_1,m)}$ and
$I^{(h)}_U$.

For $S\subset\{1,\ldots,n\}$ and $h=n-\#S$, there is a natural map
$\varphi:Y_{U,S}^{0}\ra I^{(h)}_U$ which is defined by sending the
chain of isogenies $\CC$ to its \'etale quotient. By Lemma 4.1 of
\cite{ty} this map is finite and bijective on geometric points. By
Corollary 4.2 of {\it op.\ cit.}\ we have
\begin{align*}
  H^i_c(Y^0_{U,S}\times_{k(w)}\overline{k(w)},\CL_\xi)&\isoto
  H^i_c(I^{(h)}_U\times_{k(w)}\overline{k(w)},\CL_\xi)\\ &\isoto
  H^i_c(I^{(h)}_{U^p,(1,m)}\times_{k(w)}\overline{k(w)},\CL_\xi)^{B_h(k(w))}
\end{align*}
for each $i\in\Z_{\geq 0}$ and these isomorphisms are compatible with
the action of $G(\A^{p,\infty})$ as $U^p$ varies.

If $l\neq p$, set
\[ H^i_c(I^{(h)}_{\Iw(m)},\CL_\xi) =\varinjlim_{U^p}H^i_c(I^{(h)}_U\times_{k(w)}\overline{k(w)},\CL_\xi).\]
This is an admissible $\Qlbar[G(\A^{\infty,p})\times
\Frob_w^\Z]$-module. Define
\begin{align*}
  [H(Y_{\Iw(m),S},\CL_\xi)] & = \sum_i (-1)^{n-\#S-i}H^i(Y_{\Iw(m),S},\CL_\xi)\\
  [H_c(Y^0_{\Iw(m),S},\CL_\xi)] & = \sum_i (-1)^{n-\#S-i} H^i_c(Y^0_{\Iw(m),S},\CL_\xi)\\ 
    [H_c(I^{(h)}_{\Iw(m)},\CL_\xi)] & = \sum_i (-1)^{h-i}H^i_c(I^{(h)}_{\Iw(m)},\CL_\xi)
\end{align*}
in $\Groth_l(G(\A^{\infty,p})\times\Frob_w^\Z)$.

There is, up to isomorphism, a unique one-dimensional compatible formal \linebreak Barsotti-Tate
$\CO_{F,w}$-module $\Sigma_{F_w,n-h}$ over $\overline{k(w)}$ of $\CO_{F,w}$-height
$n-h$. We have
$\End_{\CO_{F,w}}(\Sigma_{F_w,n-h})\otimes_\Z\Q\cong D_{F_w,n-h}$, the
division algebra with centre $F_w$ and Hasse invariant $1/(n-h)$. For $m_1\in \Z_{\geq 0}$, let
$\Ig^{(h)}_{U^p,(m_1,m)}/(\barX^{(h)}_{U_0}\times_{k(w)}\overline{k(w)})$
denote the moduli space of $\CO_{F,w}$-equivariant isomorphisms
\begin{align*}
 j^{\et} : (w^{-m_1}\CO_{F,w}/\CO_{F,w})^h_{\barX^{(h)}_{U_0}\times_{k(w)}\overline{k(w)}}& \isoto
 \CG^{\et}[w^{m_1}] \\
 j^0 : (\Sigma_{F_w,n-h}[w^{m_1}])_{\barX^{(h)}_{U_0}\times_{k(w)}\overline{k(w)}} & \isoto
 \CG^0[w^{m_1}]
\end{align*}
 that extend \'etale locally to any level $m_1'>m_1$. (In the notation of \cite{shin}, for each $0\leq h\leq n-1$, there is a unique
 $b\in B(G_{\Q_p},-\mu)$ corresponding to $h$ (see displayed equation (5.3) of {\it op.\ cit.}). If $m=(m_1,\ldots,m_1)$, then
 $\Ig^{(h)}_{U^p,(m_1,m)}$ is denoted $\Ig_{b,U^p,m_1}$ in \cite{shin}
 (see Section 5.2 of {\it op.\ cit.} and Section 4 of \cite{mant}). We have simply extended the
 definition to `non-parallel' $(m_1,m)\in\Z_{\geq
   0}\times\Z^{r-1}_{\geq 0}$. We also note that the notation
 $\Ig^{(h)}$ is used in place of $\Ig_b$ in Section 7.3 of \cite{shin}.)  If $l\neq p$ and $\xi$ is an
 irreducible algebraic representation of $G$ over $\Qlbar$, then $\xi$
 gives rise to a lisse sheaf $\CL_\xi$ on each
 $\Ig^{(h)}_{U^p,(m_1,m)}$. Let 
\[
H^i_c(\Ig^{(h)},\CL_\xi)=\varinjlim_{U^p,m_1,m}H^i_c(\Ig^{(h)}_{U^p,(m_1,m)},\CL_\xi).\]
This is an admissible $\Qlbar[G(\A^{\infty,p})\times
J^{(h)}(\Q_p)]$-module where
\[ J^{(h)}(\Q_p)=\Q_p^\times \times (
D_{F_w,n-h}^\times\times GL_h(F_w))\times\prod_{i=2}^r GL_n(F_{w_i})\]
(see Section 5 of \cite{shin}, where $J^{(h)}(\Q_p)$ is denoted
$J_b(\Q_p)$, with $b$ being the element of $B(\Q_p,-\mu)$ corresponding to $h$; in Section
7.3 of {\it op.\ cit.}, $J_b$ is denoted $J^{(h)}$). We have 
\[ H^i_c(\Ig^{(h)},\CL_\xi)^{\Z_p^\times \times (
  \CO_{D_{F_w,n-h}}^\times\times \Iw_{h,w})\times
  U_p^w(m)}\cong H^i_c(I^{(h)}_{\Iw(m)},\CL_\xi),\]
where the latter is regarded as an admissible
$\Qlbar[G(\A^{p,\infty})]$-module. Moreover, the action of $\Frob_w$ on the
right hand side corresponds to the action of 
\[(1,p^{-[k(w):\F_p]},\varpi_{D_{F_w,n-h}}^{-1},1,1)\in
G(\A^{p,\infty})\times \Q_p^\times \times D_{F_w,n-h}^\times \times
GL_h(F_w)\times \prod_{i=2}^r GL_n(F_{w_i})\]
on the left hand side, where $\varpi_{D_{F_w,n-h}}$ is any uniformizer
in $D_{F_w,n-h}$. We let
\[ [H_c(\Ig^{(h)},\CL_\xi)] = \sum_i
(-1)^{h-i}H^i_c(\Ig^{(h)},\CL_\xi) \]
in $\Groth_l(G(\A^{\infty,p})\times J^{(h)}(\Q_p))$.
As on page 489 of \cite{ty}, we have 
\[ [H(Y_{\Iw(m),S},\CL_\xi)]=\sum_{T\supset S}(-1)^{(n-\#S)-(n-\#T)}[H_c(I^{(n-\#T)}_{\Iw(m)},\CL_\xi)].\]
As there are $\begin{pmatrix}  n-\#S \\ h \end{pmatrix}$ subsets
$T\supset S$ with $n-\#T=h$, we deduce the following:

\begin{lem}
\label{lem:cohom-Y-and-Ig}
Suppose $l\neq p$ and $S\subset\{1,\ldots,n\}$. Then we have an equality
\[ \begin{array}{l} [H(Y_{\Iw(m),S},\CL_\xi)]= \\ \sum_{h=0}^{n-\#S}(-1)^{n-\#S-h}
\begin{pmatrix}  n-\#S \\ h \end{pmatrix}
[H_c(\Ig^{(h)},\CL_\xi)]^{\Z_p^\times \times (
  \CO_{D_{F_w,n-h}}^\times\times \Iw_{h,w})\times
  U_p^w(m)} \end{array} \]
in $\Groth_l(G(\A^{p,\infty})\times \Frob_w^\Z)$.
\end{lem}

\section{Computing the cohomology of $Y_{U,S}$}
\label{subsec:comp-cohom-Y}

In this section we deduce analogues of Proposition 4.4 of \cite{ty}.

\subsection{The stable case}
\label{subsubsec:cohom-stable-case}

Let $\Pi^1$ be a RACSDC automorphic representation of the group
$GL_n(\A_F)$. Suppose that $\Pi^1$ is $\Xi^1$-cohomological where
$\Xi^1$ is an irreducible algebraic representation of $R_{F/\Q}GL_n$
over $\C$. Assume that
\begin{itemize}
\item $\Ram_\Q(\Pi^1)\subset \Spl_{F/F^+,\Q}$.
\end{itemize}
By Lemma 7.2 of \cite{shin}, we can and do choose an
algebraic Hecke character $\psi : \A_E^\times/E^\times\ra\C^\times$
and an algebraic representation $\xi_\C$ of $G$ over $\C$ such that
\begin{itemize}
\item $\psi_{\Pi^1}|_{\A_E^\times}=\psi^c/\psi$;
\item If $\Xi$ is the representation of $\G_n$ over $\C$ corresponding
  to $\xi_\C$ as in Section \ref{subsec:notat-runn-assumpt}, then $\Xi^1$ is isomorphic to the
  restriction of $\Xi$ to $(R_{F/Q}GL_n)\times_\Q\C$;
\item $\xi_\C|^{-1}_{E_\infty^\times}=\psi_\infty^c$ (see below);
\item $\Ram_{\Q}(\psi)\subset \Spl_{F/F^+,\Q}$;
\item $\psi$ is unramified at $u$ (recall that $u$ is the prime of $E$
below the $w_i$).
\end{itemize}
(We note that Lemma 7.2 of \cite{shin} does not guarantee that $\psi$
be unramified at $u$, but the fact that this can be achieved follows
from the proof of Lemma VI.2.10 of \cite{ht}.)  In the third bullet
point, we consider $E_\infty^\times$ embedded in $G(\R)\subset
\R^\times \times GL(V\otimes_{\Q}\R)$ via the map $z \mapsto
(zz^c,z)$. It then follows from the third bullet point that
$R_{l,\imath}(\psi)$ is pure of weight $m_\xi-2t_\xi$. Set
\[ \Pi:=\psi\otimes \Pi^1. \]
Then $\Pi$ is a $\Xi$-cohomological automorphic representation of
  $\G_n(\A)\cong GL_1(\A_E)\times GL_n(\A_F)$. 
Note that $\Pi^1$, $\psi$ and $\Pi$ satisfy the assumptions of Section
\ref{subsubsec:stable-case}. Let $\xi = \imath^{-1}\xi_\C$, an irreducible algebraic
representation of $G$ over $\Qlbar$. Let
$\pi_p\in\Irr_l(G(\Q_p))$ be such that $\BC_p(\imath\pi_p)\cong\Pi_p$
(note that $\pi_p$ is unique as $p$ splits in $E$).

The next result follows from Proposition 7.14 of \cite{shin}.

\begin{prop}
  \label{prop: cohomology vanishing stable case}
  Suppose that $l\neq p$ and $\pi_p^{\Iw(m)}\neq(0)$. Let
  $T\supset\{\infty\}$ be a finite set of places of $\Q$ with
  $\Ram_{F/\Q}\cup \Ram_{\Q}(\Pi)\cup\{p\}\subset T_{\fin}\subset
  \Spl_{F/F^+,\Q}$.  Then for every $S\subset \{1,\ldots,n\}$, we have
\[ H^j(Y_{\Iw(m),S},\CL_\xi)^{K^T}\{\Pi^T\} =(0) \]
for $j\neq n-\#S$.
\end{prop}

The following corollary can be proved in the same way as Corollary 4.5 of \cite{ty}.

\begin{corollary}
 \label{cor:cohom vanishing stable l=p}
  Suppose that $l=p$ and $\sigma :W_0\into \Qlbar$ over
  $\Z_p=\Z_l$. Let $T\supset\{\infty\}$ be a finite set of places of
  $\Q$ with $\Ram_{F/\Q}\cup \Ram_{\Q}(\Pi)\cup\{l\}\subset
  T_{\fin}\subset \Spl_{F/F^+,\Q}$. If $\pi_l^{\Iw(m)}\neq(0)$, then
  for every $S\subset\{1,\ldots,n\}$, we have
\[ a_\xi (H^{j+m_\xi}(\CA_{\Iw(m),S}^{m_\xi}/W_0)\otimes_{W_0,\sigma}\Qlbar)^{K^T}\{\Pi^T\} =(0) \]
for $j\neq n-\#S$.
\end{corollary}

In the next result we place no restriction on the primes $l$ and $p$.

\begin{corollary}
  \label{cor:purity-stable} 
  If $\pi_p^{\Iw(m)}\neq 0$, then
  $\WD(\tR^{n-1}_{\xi,l}(\Pi)|_{G_{F_w}})$ is pure of weight
  $m_\xi-2t_\xi+n-1$ and $\WD(R_{l,\imath}(\Pi^1)|_{G_{F_w}})$ is pure of weight $n-1$.
\end{corollary}

\begin{proof}
  Let
  $T=\{\infty\}\cup\Ram_{F/\Q}\cup\Ram_{\Q}(\Pi)\cup\{p\}$
  and let $D=\dim \pi_p^{\Iw(m)}$. Let $T'=T_{\fin}-\{p\}$.
  By Theorem \ref{thm:cohomology-stable-case}, we have an isomorphism
  of $\Qlbar[\G_n(\A_{T'})\times G_F]$-modules
  \[ \BC_{T'}(H^{n-1}(X_{\Iw(m)},\CL_\xi)^{K^T}\{\Pi^T\})
  \cong  ((\imath^{-1}\Pi_{T'})\otimes
  \tR^{n-1}_{\xi,l}(\Pi))^{\oplus D}. \]
  By proposition \ref{prop:RZ-spec-seq}, there is a spectral sequence
 \[ E^{i,j}_1(\Iw(m),\xi)^{K^T}\{\Pi^T\}\implies
 \WD(H^{n-1}(X_{\Iw(m)},\CL_\xi)|_{G_{F_w}})^{K^T}\{\Pi^T\}.\]
 Using Proposition \ref{prop: cohomology vanishing stable case} (when $l\neq
 p$) and Corollary
 \ref{cor:cohom vanishing stable l=p} (when $l=p$) we see that
 $E^{i,j}_1(\Iw(m),\xi)^{K^T}\{\Pi^T\}=(0)$ unless $i+j=n-1$, and thus
 the spectral sequence degenerates at $E_1$. Let $\pi_{T'}$
 denote the unique element of $\Irr_l(G(\A_{T'}))$ with $\BC_{T'}(\pi_{T'})=\imath^{-1}\Pi_{T'}$.
 Then, for $i+j=n-1$,
 $E_1^{i,j}(\Iw(m),\xi)^{K^T}\{\Pi^T\}$ is
 of the form $\pi_{T'}\otimes R_j$ where $R_j$ is a finite dimensional
 $\Qlbar[\Frob_w^\Z]$-module which is pure of weight
 $j+m_\xi-2t_\xi$ (and possibly zero). The first statement now
 follows from this and the description of the monodromy operator $N$
 in Proposition \ref{prop:RZ-spec-seq}. The second statement follows
 from the first statement together with Theorem
 \ref{thm:cohomology-stable-case} and Lemma 1.7 of \cite{ty}.
\end{proof}

\subsection{The endoscopic case}
\label{subsubsec:comp-cohom-Y-endo}

Suppose we are in the following situation:
\begin{itemize}
\item $\Pi_1$ is a RACSDC automorphic representation of $GL_{n-1}(\A_F)$;
\item $\Ram_\Q(\Pi_1)\subset \Spl_{F/F^+,\Q}$;
\item $\Pi_1$ is cohomological for an irreducible algebraic
  representation $\Xi_1$ of the group $R_{F/\Q}(GL_{n-1})$ over $\C$;
\item $\Pi_1$ has Shin-regular weight.
\end{itemize}

\begin{lem}
\label{lem:aux data in endoscopic case}
We can find
\begin{itemize}
\item a continuous algebraic character
  $\Pi_2:\A_F^\times/F^\times\ra\C^\times$ with $\Pi_2^{-1}=\Pi_2\circ
  c$;
\item a continuous algebraic character $\psi :
  \A_E^\times/E^\times\ra\C^\times$;
\item a continuous character $\varpi :
  \A_E^\times/E^\times\ra\C^\times$; and
\item an irreducible algebraic representation $\xi_\C$ of $G$ over
  $\C$
\end{itemize} 
such that if we set
\begin{align*}
\Pi_{M,1}&:=\Pi_1\otimes(\varpi\circ\norm_{F/E}\circ\det)\\
\Pi_{M,2}&:=\Pi_2\\
\Pi^1&:=\nind_{GL_{n-1}\times GL_1}^{GL_n}(\Pi_{M,1}\otimes\Pi_{M,2}) 
\end{align*}
and let
\begin{itemize}
\item $\Xi$ be the irreducible algebraic representation of $\G_n$
over $\C$ which corresponds to $\xi_\C$ as in Section \ref{subsec:notat-runn-assumpt};
\item $\Xi^1:=\Xi|_{R_{F/Q}(GL_n)\times_\Q\C}$,
\end{itemize}
then
\begin{itemize}
\item $\Ram_\Q(\Pi_2)\subset \Spl_{F/F^+,\Q}$;
\item $\Pi_2$ is unramified at $u$;
\item $\Ram_\Q(\psi)\subset \Spl_{F/F^+,\Q}$;
\item $\xi_\C|_{E_{\infty}^\times}^{-1}=\psi_{\infty}^c$;
\item $\psi$ is unramified at $u$;
\item $\Ram_\Q(\varpi)\subset \Spl_{F/F^+,\Q}$;
\item $\varpi|_{\A^\times}$ factors through
  $\A^\times/\Q^\times\R^\times_{>0}$ and equals the composite of
  $\Art_{\Q}$ with the surjective character
  $G_{\Q}^{\ab}\onto\Gal(E/\Q)\iso\{\pm 1\}$ (note that this implies
  $\varpi^{-1}=\varpi\circ c$);
\item $\varpi$ is unramified at $u$;
\item $\Pi^1$ is cohomological for $\Xi^1$ (note that $\Pi^1$ is
  irreducible, as $\Pi_{M,1}$ and $\Pi_{M,2}$ are unitary, and also that
  $(\Pi^1)^\vee\cong\Pi^1\circ c$);
\item $\psi_{\Pi^1}|_{\A_E^\times}=\psi^c/\psi$ (recall that $\psi_{\Pi^1}$
  denotes the central character of $\Pi^1$).
\end{itemize}
Moreover, if we apply Theorem \ref{thm:cohomology-endo-case} to $\Pi$,
then alternative (2)(a) holds. In other words, the integer $e_2(\Pi,G)$
equals $1$ and if $\xi =
\imath^{-1}\xi_\C$, then
\[ \tR_{\xi,l}^{n-1}(\Pi)^{\semis}\cong R_{l,\imath}(\Pi_1)^{G_G}\otimes
R_{l,\imath}(\psi\varpi |\cdot|^{1/2})|_{G_F}.\]
\end{lem}

\begin{proof}
  This follows by combining Lemmas 7.1, 7.2 and 7.3 of
  \cite{shin}. (More precisely, we first choose $\varpi$ using Lemma
  7.1. The extra condition that $\varpi$ be unramified at $u$ is
  easily achieved -- in the proof of Lemma 7.1 we add the primes $u$
  and $u^c$ to the set $R$ and insist that $\varpi^0$ takes value $1$
  on $p\in E_u^\times$ and on $p\in E_{u^c}^\times$. We then make two
  candidate choices $\chi$ and $\chi'$ for $\Pi_2$ with
  $\chi^{-1}=\chi\circ c$, $\Ram_{\Q}(\chi)\subset \Spl_{F/F^+}(\Q)$
  and $\chi$ unramified at $u$ and with $\chi'$ having the same
  properties. In addition, we assume that the infinity type of $\chi$
  and $\chi'$ are as prescribed in the paragraph before Lemma 7.2 of
  \cite{shin}. (The fact that we can find such characters follows for
  instance from Lemma 2.2 of \cite{hsbt}.) Lemma 7.2 of \cite{shin}
  then tells us that we can choose pairs $(\psi,\xi_\C)$ and
  $(\psi',\xi'_\C)$ corresponding to the choice of $\Pi_2=\chi$ or
  $\Pi_2=\chi'$ and satisfying all the required properties except the
  requirement that $\psi$ and $\psi'$ be unramified at $u$ and the
  requirement that the integer $e_2(\Pi,G)$ equal $1$. However, the
  proof of Lemma VI.2.10 of \cite{ht} shows that we may choose $\psi$
  and $\psi'$ to be unramified at $u$ and Lemma 7.3 of \cite{shin} shows
  that for one of the choices $(\chi,\psi,\xi_\C)$ or
  $(\chi',\psi',\xi'_\C)$, the corresponding integer $e_2(\Pi,G)$ equals
  1.)
\end{proof}

Choose $\Pi_2$, $\psi$, $\varpi$ and $\xi_{\C}$ as in the above lemma and
keep all additional notation introduced there. Let
\[ \Pi = \psi \otimes \Pi^1, \] an automorphic representation of
$\G_n(\A)$. Let $\pi_p\in\Irr_l(G(\Q_p))$ be the unique representation with
$\BC(\imath\pi_p)\cong\Pi_p$. Write
$\pi_p=\pi_p^0\otimes\pi_w\otimes(\otimes_{i=2}^r\pi_{w_i})$
corresponding to the decomposition $G(\Q_p)\cong \Q_p^\times \times
GL_n(F_w)\times \prod_{i=2}^r GL_n(F_{w_i})$. Let $\pi_{M,w,1}=\imath^{-1}\Pi_{M,1,w}\in
\Irr_l(GL_{n-1}(F_w))$ and $\pi_{M,w,2}=\imath^{-1}\Pi_{M,2,w}\in\Irr_l(GL_1(F_w))$.

For $a,b\in \Z_{\geq 0}$, let 
\[ \nRed^{a,b} :
\Groth_l(GL_{a+b}(F_w))\ra\Groth_l(D_{F_w,a}^\times\times GL_b(F_w)) \]
denote the composition
\[ \begin{array}{rcl}
\Groth_l(GL_{a+b}(F_w)) & \xrightarrow{J_{N_b^\op}} & \Groth_l(GL_{a}(F_w)\times
GL_{b}(F_w))\\ & \xrightarrow{\LJ_{a}\otimes \id} & \Groth_l(D_{F_w,a}^\times\times GL_b(F_w))
\end{array} \]
where $N_b^{\op}$ is unipotent radical of the parabolic subgroup of $GL_m$ consisting
of block lower triangular matrices with an $(a\times a)$-block in the
upper left corner and a $(b\times b)$-block in the lower right corner;
\[ J_{N_b^{\op}}:\Groth_l(GL_m(F_w))\ra\Groth_l(GL_a(F_w)\times GL_b(F_w)) \]
 is the normalized
Jacquet module functor; and 
\[ \LJ_a : \Groth_l(GL_a(F_w))\ra \Groth_l(D_{F_w,a}^\times) \]
  is the map denoted $\mathbf{LJ}_1$ in Proposition 3.2 of \cite{badu}.
(See Section 2.4 of \cite{shin}.)

Let
\[ \overline{\delta}_{P_h}^{1/2} : J^{(h)}(\Q_p)\ra
\C^\times \]
denote the character which sends $(g_{p,0},(d,g),g_i)\in \Q_p^\times
\times (D_{F_w,n-h}^\times\times GL_h(F_w))\times\prod_{i=2}^r
GL_n(F_{w_i})$ to $|\det(d)^h\det(g)^{-(n-h)}|_{F_w}^{1/2}$.

\begin{thm}
\label{thm:cohom of igusa variety}
Suppose $l\neq p$. Let
  $T\supset\{\infty\}$ be a finite set of places of $\Q$ with
  $\Ram_{F/\Q}\cup\Ram_{\Q}(\Pi)\cup\Ram_\Q(\varpi)\cup\{p\}\subset T_{\fin}\subset
  \Spl_{F/F^+,\Q}$.
 Then
\[ [H_c(\Ig^{(0)},\CL_\xi)][\Pi^T] = (0) \]
while for $1\leq h \leq n-1$, we have an equality
{\small \begin{align*}
  &\BC^p([H_c(\Ig^{(h)},\CL_\xi)][\Pi^T])=
  C_G[\imath^{-1}\Pi^{\infty,p}] \times \\
 &\left[\left(\pi_{p,0}\otimes\nind_{GL_{h-1,1}(F_w)}^{GL_h(F_w)}(\nRed^{n-h,h-1}(\pi_{M,w,1})\otimes\pi_{M,w,2})\otimes
  (\otimes_{i=2}^r \pi_{w_i})\right)\otimes \imath^{-1}\overline{\delta}^{1/2}_{P_h}\right]
\end{align*} }
in $\Groth_l(\G_n(\A^{\infty,p})\times J^{(h)}(\Q_p))$.
\end{thm}

\begin{proof}  
  The result is essentially a rewording of part (ii) of Theorem 6.1 of
  \cite{shin}. We freely make use of the notation of {\it op.\ cit.}\
  for the rest of this proof. Let $0\leq h\leq n-1$ and let $b\in
  B(\Q_p,-\mu)$ correspond to
  $h$ (in the sense explained above). The constant $e_1$ which appears
  in the statement of Theorem 6.1 of \cite{shin} is equal to $(-1)^{n-1}=1$ by Corollary
  6.5(ii) of {\it op.\ cit}. Lemma \ref{lem:aux data in endoscopic
    case} above, and the choices made after it, guarantee that the constant $e_2$ also equals
  $1$. Applying Theorem 6.1 of
  \cite{shin}, we obtain
\[ \begin{array}{l} \BC^p([H_c(\Ig^{(h)},\CL_\xi)][\Pi^T])= \\
  C_G(-1)^h[\imath^{-1}\Pi^{\infty,p}] \times \left[
    \frac{1}{2}\left(\Red^b_n(\pi_p)+\Red^b_{n-1,1}(\pi_{H,p})\right)\right]. \end{array} \]
(We remark that our definition of
$[H_c(\Ig^{(h)},\CL_\xi)]$ differs from Shin's definition of
$H_c(\Ig_b,\CL_\xi)$ by a factor of $(-1)^h$.) We have $\Red^b_{n-1,1}(\pi_{H,p})=\nRed^b_{n-1,1}(\pi_{H,p})\otimes
\imath^{-1}\overline{\delta}^{1/2}_{P_h}$ (see Section 5.5 of
\cite{shin}). By Lemma 5.9 of \cite{shin} and the discussion
immediately preceding it, we have
\[ \nRed^b_n(\pi_p)+\nRed^b_{n-1,1}(\pi_{H,p}) =e_p(J_b)\pi_{p,0}\otimes
2X_1(h,\pi_{H,p})\otimes(\otimes_{i=2}^r \pi_{w_i}) \]
where
\[ X_1(h,\pi_{H,p}) =
\begin{cases} 0 & (h=0) \\
\nind_{GL_{h-1,1}(F_w)}^{GL_h(F_w)}(\nRed^{n-h,h-1}(\pi_{M,w,1})\otimes\pi_{M,w,2})
& (h\neq 0)
\end{cases}.\]
Moreover, $e_p(J_b)=(-1)^{n-h-1}=(-1)^h$ (see Case 1 in Section 5.5 of \cite{shin}).
The result follows.
\end{proof}

In Remark 7.16 of \cite{shin}, Shin indicates that the following
result can proved in the same way as Proposition 7.14 of {\it op.\
  cit.}. We give a self-contained proof here for the benefit of the reader.

\begin{prop}
\label{prop: cohomology vanishing endo case}
Suppose $l\neq p$ and $\pi_p^{\Iw(m)}\neq (0)$. Let
$T\supset\{\infty\}$ be a finite set of places of $\Q$ with
$\Ram_{F/\Q}\cup\Ram_{\Q}(\Pi)\cup\Ram_\Q(\varpi)\cup\{p\}\subset T_{\fin}\subset
\Spl_{F/F^+,\Q}$.  Then for every $S\subset \{1,\ldots,n\}$ we have
\[ H^j(Y_{\Iw(m),S},\CL_\xi)^{K^T}\{\Pi^T\}=(0) \]
for $j \neq n-\# S$.
\end{prop}

\begin{proof}
Let $D=C_G(\dim
(\otimes_{i=2}^r \pi_{w_i})^{U^w_p(m)})$. We deduce from Theorem
\ref{thm:cohom of igusa variety} and Lemma \ref{lem:cohom-Y-and-Ig} that
{\small \begin{align*}
&\BC^p [H(Y_{\Iw(m),S},\CL_\xi)][\Pi^T]=D[\imath^{-1}\Pi^{\infty,p}]\times \sum_{h=1}^{n-\# S}(-1)^{n-\# S-h}
\begin{pmatrix} n-\# S\\ h\end{pmatrix}\\
&\left[\pi_{p,0}^{\Z_p^\times}\otimes
\left(\nind_{GL_{h-1,1}(F_w)}^{GL_h(F_w)}(\nRed^{n-h,h-1}(\pi_{M,w,1})\otimes\pi_{M,w,2})\otimes\imath^{-1}\overline{\delta}_{P_h}^{1/2}\right)^{\CO_{D_{F_w,n-h}}^\times\times\Iw_{h,w}}
\right]
\end{align*}}
in $\Groth_l(\G_n(\A^{p,\infty})\times\Frob_w^\Z)$ (recall that
$\Frob_w$ acts via
$(p^{-[k(w):\F_p]},(\varpi_{D_{F_w,n-h}}^{-1},1))$ in $\Q_p^\times\times
(D_{F_w,n-h}^\times\times GL_h(F_w))$).
Since $\pi_p^{\Iw(m)}\neq (0)$, we can write
$\pi_{M,w,1}=\Sp_{s_1}(\pi_1)\boxplus\dots\boxplus\Sp_{s_t}(\pi_t)$
where each $\pi_j$ is an unramified character
$\pi_j:F_w^\times\ra\Qlbar^\times$. As $\Pi_{M,1}$ is generic, we
know that
$\pi_{M,w,1}=\nind_{P(F_w)}^{GL_{n-1}(F_w)}(\Sp_{s_1}(\pi_1)\otimes\dots\otimes\Sp_{s_t}(\pi_t))$
where $P\subset GL_{n-1}$ is an appropriate parabolic subgroup.

Using Lemma I.3.9 of \cite{ht} and Th\'eor\`eme 3.1 and Proposition
3.2 of \cite{badu}, we see that for $1\leq h \leq n-\# S$,
\[ \begin{array}{l}
\nind_{GL_{h-1,1}(F_w)}^{GL_h(F_w)}(\nRed^{n-h,h-1}(\pi_{M,w,1})\otimes
\pi_{M,w,2})= \\ \sum_i
[(\pi_i\cdot \imath^{-1}|\cdot|_{F_w}^{(n-h-1)/2})\circ \det]\times 
[\nind_{P'(F_w)}^{GL_h(F_w)}(\Sp_{s_i-(n-h)}(\pi_i|\cdot|_{F_w}^{n-h})\otimes \\
\multicolumn{1}{r}{(\otimes_{j\neq
  i}\Sp_{s_j}(\pi_j))\otimes\pi_{M,w,2})]}
\end{array} \]
in $\Groth_l(D_{F_w,n-h}^\times\times GL_h(F_w))$, where the sum is
over all $i$ such that $s_i\geq n-h$, and $P'\subset
GL_h$ is an appropriate parabolic subgroup. (As pointed out in
\cite{ty}, there is a typo in Lemma I.3.9 of \cite{ht} --- `positive integers
$h_1,\ldots,h_t$' should be replaced by `non-negative integers $h_1,\ldots,h_t$'.) Moreover,
\begin{align*} 
&\dim\left( \nind_{P'(F_w)}^{GL_h(F_w)}(\Sp_{s_i-(n-h)}(\pi_i|\cdot|_{F_w}^{n-h})\otimes(\otimes_{j\neq
  i}\Sp_{s_j}(\pi_j))\otimes\pi_{M,w,2})\right)^{\Iw_{h,w}}\\
&=\frac{h!}{(s_i-(n-h))!\prod_{j\neq i}s_j!}
\end{align*}
(see page 490 of \cite{ty}). From this we deduce that
\begin{align*} 
&\BC^p [H(Y_{\Iw(m),S},\CL_\xi)][\Pi^T]=D[\imath^{-1}\Pi^{\infty,p}]\times\\ 
&\sum_{h=1}^{n-\# S}(-1)^{n-\# S-h}
\begin{pmatrix} n-\# S\\ h\end{pmatrix}
\sum_{i:s_i\geq n-h}\frac{h!}{(s_i-(n-h))!\prod_{j\neq i}s_j!}[V_i],
\end{align*}
where $V_i=\rec(\pi_i^{-1}\cdot\imath^{-1}|\cdot|_{F_w}^{(1-n)/2}(\pi_{p,0}\circ
\norm_{F_w/E_u})^{-1})$. As on page 490 of \cite{ty}, it follows that
\[ \BC^p
[H(Y_{\Iw(m),S},\CL_\xi)][\Pi^T]=D[\imath^{-1}\Pi^{\infty,p}]\times
\sum_{i:s_i=\# S} \frac{(n-\#S)!}{\prod_{j\neq i}s_j!}[V_i].\]
As $\Pi_{1,w}$ is unitary and tempered (by Corollary 1.3 of \cite{shin}) and
$\rec(\pi_{p,0})\cong \imath^{-1}\rec(\psi_u)$  is strictly pure of weight $2t_\xi-m_\xi$ (since
$\xi_\C|_{E_{\infty}^\times}^{-1}=\psi_\infty^c$), we see that
$\rec(\pi_{M,w,1}^\vee\otimes\imath^{-1}|\cdot|_{F_w}^{(1-n)/2}(\pi_{p,0}\circ
\norm_{F_w/E_u})^{-1})$
is pure of weight $m_\xi-2t_\xi+n-1$. If $s_i=\#S$, it follows that
$V_i$ is strictly pure of weight
$m_\xi-2t_\xi+n-1-(\#S-1)=m_\xi-2t_\xi+n-\#S$. The Weil conjectures
now imply that if $j\neq n-\#S$ then
\[ [H^j(Y_{\Iw(m),S},\CL_\xi)][\Pi^T]= \\ a_\xi
[H^{m_\xi+j}(\CA_{\Iw(m),S},\Qlbar(t_\xi))][\Pi^T] =(0) \]
 in $\Groth_l(G(\A^{p,\infty})\times\Frob_w^\Z)$. Since
\[ [H^j(Y_{\Iw(m),S},\CL_\xi)^{K^T}\{\Pi^T\}]=([H^j(Y_{\Iw(m),S},\CL_\xi)][\Pi^T])^{K^T}
\]
in $\Groth_l(G(\A_{T_{\fin}-\{p\}})\times\Frob_w^\Z)$, the result follows.
\end{proof}

The proof of Corollary 4.5 of \cite{ty} allows us to deduce
the following.

\begin{corollary}
\label{cor:cohom vanishing endo l=p}
  Suppose that $l=p$ and $\sigma :W_0\into \Qlbar$ over
  $\Z_p=\Z_l$. Let
  $T\supset\{\infty\}$ be a finite set of places of $\Q$ with
  $\Ram_{F/\Q}\cup \Ram_{\Q}(\Pi)\cup\{p\}\subset T_{\fin}\subset
  \Spl_{F/F^+,\Q}$. If $\pi_p^{\Iw(m)}\neq(0)$, then for every $S\subset\{1,\ldots,n\}$, we have
\[ a_\xi (H^{j+m_\xi}(\CA_{\Iw(m),S}^{m_\xi}/W_0)\otimes_{W_0,\sigma}\Qlbar)^{K^T}\{\Pi^T\} =(0) \]
for $j\neq n-\#S$.
\end{corollary}

The next corollary follows from the previous two results combined with Theorem
\ref{thm:cohomology-endo-case} and the proof of Corollary
\ref{cor:purity-stable}.

\begin{corollary}
\label{cor:purity-endo}
  If $\pi_p^{\Iw(m)}\neq 0$, then
  $\WD(\tR^{n-1}_{\xi,l}(\Pi)|_{G_{F_w}})$ is pure of weight
  $m_\xi-2t_\xi+n-1$ and $\WD(R_{l,\imath}(\Pi_1)|_{G_{F_w}})$ is pure
  of weight $n-2$. 
\end{corollary}

\section{Proof of Theorem \ref{thm:semistable Shin-regular case}}
\label{sec:proof-theorem-semistable-slightly-regular}

We now complete the proof of Theorem \ref{thm:semistable Shin-regular
  case}. Suppose that $v|l$ is a place of $L$ with
$\Pi_{v}^{\Iw_{m,v}}\neq\{0\}$. Choose a finite CM soluble Galois
extension $F/L$ such that
\begin{itemize}
\item $[F^+:\Q]$ is even;
\item $F=EF^+$ where $E$ is a quadratic imaginary field in which $l$ splits;
\item $F$ splits completely above $v$;
\item $\BC_{F/L}(\Pi)$ is cuspidal;
\item $\Ram_{F/\Q}\cup \Ram_{\Q}(\BC_{F/L}(\Pi))\subset \Spl_{\Q,F/F^+}$.
\end{itemize}
(See the argument in the penultimate paragraph of \cite{shin}.) Let
$l$ and $\imath:\Qlbar \iso \C$ as given to us by the statement of
Theorem \ref{thm:semistable Shin-regular case}. Choose another prime
$l'\neq l$ and $\imath':\Qbar_{l'}\iso\C$.  Recall that in Section
\ref{subsec:notat-runn-assumpt} we introduced notation that was then
in force from Section \ref{subsec:shimura-varieties} to Section
\ref{subsec:comp-cohom-Y}. We will shortly apply the results of these
sections in two scenarios -- one where the pair $(l,\imath)$ of
Section \ref{subsec:notat-runn-assumpt} is equal to the pair
$(l,\imath)$ of the statement of Theorem \ref{thm:semistable
  Shin-regular case} and one where the $(l,\imath)$ of Section
\ref{subsec:notat-runn-assumpt} is equal to the pair $(l',\imath')$
chosen above. The rest of the notation we fix as follows: we take
$E,F$ and $F^+$ as chosen above. We take $p=l$ and let $w$ be a prime
of $F$ lying above the prime $v$ of $L$. This determines $u$ and
$w_1,\ldots,w_r$. We choose some $\tau:F\into\C$ and
$\imath_p:\Qpbar\iso\C$ such that $\imath_p^{-1}\circ \tau$ induces
$w$. We take $n=m$ if $m$ is odd and $n=m+1$ otherwise. Finally we
choose a set of data $(V,\seq{\cdot,\cdot},h)$ satisfying the
assumptions of Section \ref{subsec:notat-runn-assumpt}.

Suppose first of all that $m$ is odd. Denote $\BC_{F/L}(\Pi_L)$ by
$\Pi^1$. We choose $\psi$ and $\xi_\C$ as in Section
\ref{subsubsec:cohom-stable-case} and set $\tPi=\psi\otimes\Pi^1$,
$\xi=\imath^{-1}\xi_\C$ and $\xi'=(\imath')^{-1}\xi_\C$. Define
$\tR^{n-1}_{\xi,l}(\tPi)$ and $\tR^{n-1}_{\xi',l'}(\tPi)$ as in Section
\ref{subsubsec:stable-case}. Then $\tR^{n-1}_{\xi',l'}(\tPi)^{\semis}
\cong R_{l',\imath'}(\Pi^1)^{C_G}\otimes R_{l',\imath'}(\psi)|_{G_F}$
by Theorem \ref{thm:cohomology-stable-case} and hence
\[
\imath'\WD(\tR^{n-1}_{\xi',l'}(\tPi)^{\semis}|_{G_{F_w}})^{\Fsemis}
\cong \rec((\Pi^1_w)^{\vee}\otimes |\cdot|^{(1-n)/2}\circ \det)^{C_G}
\otimes \rec(\psi_u^{-1}\circ \norm_{F_w/E_u}) \] by Theorem 1.2 of
\cite{shin}. Let $T\supset\{\infty\}$ be a finite set of places of
$\Q$ with $\Ram_{F/\Q}\cup\Ram_{\Q}(\tPi)\cup\{p\}\subset
T_{\fin}\subset \Spl_{F/F^+,\Q}$ and let $T'=T_{\fin}-\{p\}$. Let
$\pi'_{T_{\fin}}$ be the unique element of
$\Irr_{l'}(G(\A_{T_{\fin}}))$ with
$\BC_{T_{\fin}}(\imath'\pi'_{T_{\fin}})=\tPi_{T_{\fin}}$. Choose
$m\in\Z^{r-1}$ and a compact open subgroup $U_{T'}\subset G(\A_{T'})$
such that $(\pi'_{T_{\fin}})^{\Iw(m)\times U_{T'}} \neq\{0\}$.  Let
$(e')^T\in \Qbar_{l'}[K^T\backslash G(\A^T)/K^T]$ be an idempotent
with $(e')^T R^{K^T} = R^{K^T}\{\tPi^T\}$ whenever $R$ is one of
$H^j(X_{\Iw(m)},\CL_{\xi'})$ or
$a_{\xi'}H^j(\CA_{\Iw(m),S}^{m_{\xi'}},\Qbar_{l'})$. Then each of
these spaces is $\pi'_{T'}$-isotypic. Let $e = \imath^{-1}\imath'
e'$. Then for each $\alpha \in W_{F_w}$, $j\geq 0$, $S\subset\{1,\ldots,n\}$ and $\sigma : W_0 \into
\Qlbar$ over $\Z_l$, we have
\[ \begin{array}{l} 
\imath'\tr(\alpha e'a_{\xi'}|H^j(\CA_{\Iw(m),S}^{m_{\xi}},\Qbar_{l'})^{K^T\times U_{T'}}) = \\
\imath \tr(\alpha ea_\xi|
(H^j(\CA_{\Iw(m),S}^{m_{\xi'}}/W_0)\otimes_{W_0,\sigma}\Qlbar)^{K^T\times U_{T'}}) \end{array}
\]
by the main results of \cite{km} and \cite{gm}. For each $j\geq 0$, we
have
\begin{align*}
  ea_\xi (H^j(\CA_{\Iw(m),S}^{m_{\xi'}}/W_0)\otimes_{W_0,\sigma}\Qlbar)^{K^T} &\subset &a_\xi
(H^j(\CA_{\Iw(m),S}^{m_{\xi'}}/W_0)\otimes_{W_0,\sigma}\Qlbar)^{K^T}\{\tPi^T\} \\
e (H^j(X_{\Iw(m)}/W_0)\otimes_{W_0,\sigma}\Qlbar)^{K^T} &\subset &
(H^j(X_{\Iw(m)}/W_0)\otimes_{W_0,\sigma}\Qlbar)^{K^T}\{\tPi^T\}.
\end{align*}
We then deduce from the previous equality of traces together with
Proposition \ref{prop: cohomology vanishing stable case}, Corollary \ref{cor:cohom vanishing stable l=p},
Proposition \ref{prop:RZ-spec-seq} and Theorem
\ref{thm:cohomology-stable-case} that the two inclusions above are
equalities (for dimension reasons) and moreover that
\[ \imath'\tr(\alpha|\WD( \tR^{n-1}_{\xi',l'}(\tPi)|_{G_{F_w}})) = \imath\tr(\alpha
| \WD(\tR^{n-1}_{\xi,l}(\tPi)|_{G_{F_w}}))\]
for each $\alpha \in W_{F_w}$ and hence
\[ \imath\WD(\tR^{n-1}_{\xi,l}(\tPi)|_{G_{F_w}})^{\semis}\cong
(\rec((\Pi^1_w)^{\vee}\otimes |\det|^{(1-n)/2})^{\semis})^{C_G}\otimes
\rec(\psi_u^{-1}\circ \norm_{F_w/E_u}).\] 
Since $\tR^{n-1}_{\xi,l}(\tPi)^{\semis} \cong
R_{l,\imath}(\Pi^1)^{C_G}\otimes R_{l,\imath}(\psi)|_{G_F}$, by Theorem
\ref{thm:cohomology-stable-case}, we see that
\[ \imath\WD(R_{l,\imath}(\Pi^1)|_{G_{F_w}})^{\semis} \cong
\rec((\Pi^1_w)^{\vee}\otimes |\det|^{(1-n)/2})^{\semis}.\] 
By Proposition \ref{thm:redn to ss}, it suffices to show
that $\WD(R_{l,\imath}(\Pi^1)|_{G_{F_w}})$ is pure and this is established in
Corollary \ref{cor:purity-stable}.  As $v$ splits completely in $F$,
we have established Theorem \ref{thm:semistable Shin-regular case}
in the case when $m$ is odd.

Now suppose that $m$ is even and denote $\BC_{F/L}(\Pi_L)$ by $\Pi_1$.
We choose $\psi$, $\xi_\C$, $\varpi$, $\Pi_2$ and $\Pi^1$ as in Lemma
\ref{lem:aux data in endoscopic case}. Set
$\tPi=\psi\otimes\Pi^1$ and $\xi=\imath^{-1}\xi_\C$ and
$\xi'=(\imath')^{-1}\xi_\C$ and define $\tR_{\xi,l}^{n-1}(\tPi)$ and
$\tR_{\xi',l'}^{n-1}(\tPi)$ as in Section
\ref{subsubsec:endoscopic-case}. The proof now proceeds exactly as
in the case where $m$ is odd except that we replace the appeals to
Theorem \ref{thm:cohomology-stable-case}, Proposition \ref{prop:
  cohomology vanishing stable case}, Corollary \ref{cor:cohom
  vanishing stable l=p} and Corollary \ref{cor:purity-stable} with
appeals to Theorem \ref{thm:cohomology-endo-case}, Proposition \ref{prop:
  cohomology vanishing endo case}, Corollary \ref{cor:cohom
  vanishing endo l=p} and Corollary \ref{cor:purity-endo} respectively.

\bibliographystyle{amsalpha}
\bibliography{barnetlambgeegeraghty}

\end{document}